\documentclass[a4paper,UKenglish,cleveref, autoref, thm-restate]{article}
\usepackage[margin=1.58in]{geometry}
\usepackage[utf8]{inputenc}
\usepackage{graphicx}
\usepackage{float}
\usepackage{dsfont,complexity,xcolor}
\usepackage{url}
\usepackage{amsmath,amsfonts,amssymb}
\usepackage{multirow}
\usepackage{wrapfig}
\usepackage{xcolor}
\usepackage{amsthm}

\usepackage{caption}
\usepackage{subcaption}

\newtheorem{clm}{Claim}
\newtheorem{theorem}{Theorem}

\newtheorem{proposition}[theorem]{Proposition}
\newtheorem{lemma}[theorem]{Lemma}
\newtheorem*{clm*}{Claim}
\newtheorem{remark}[theorem]{Remark}
\newtheorem{question}[theorem]{Question}

\usepackage{enumitem}

\usepackage{xspace}

\usepackage{soul,color}

\newcommand{\probgenarg}[1]{$k+1$ $X$-spanning $Y$-disjoint Branching\xspace}

\title{On Computing Large Temporal (Unilateral) Connected Components}

\author{Isnard Lopes Costa$^1$ \and Raul Lopes$^{2,3}$ \and Andrea Marino$^4$ \and Ana Silva$^{1,4}$}
\date{
  $^1$Universidade Federal do Ceará, Brazil.\\ \texttt{isnard.lopes@alu.ufc.br}; \texttt{anasilva@mat.ufc.br}\\%
  $^2$Université Paris-Dauphine, CNRS UMR7243, France.\\
  $^3$DIENS, École normale supérieure de Paris, CNRS, France.\\
   \texttt{raul.wayne@gmail.com}\\%
  $^4$Universitá degli Studi di Firenze, Italy.\\ \texttt{andrea.marino@unifi.it}
}

\usepackage{complexity}

\usepackage{tikz-network}
\usepackage{pgf}
\usetikzlibrary{arrows}
\usetikzlibrary{fit}
\newcommand{\vertex}{\node[circle,draw,fill=white, inner sep=2pt]}
\tikzset{my label/.style n args={2}{label={[font=\small,text=blue]#1:#2}}]}

\newcommand\myprob[3]{%
  \vspace{1mm}
  \par\noindent {\bf Problem} {#1}.\\
  {\bf Input}: {#2}\\
  {\bf Question}: {#3}\par
  \vspace{1mm}
}

\newcommand{\reach}{\textsc{R}}

\newcommand{\TCC}{\textsc{tcc}\xspace}
\newcommand{\cTCC}{\textsc{closed tcc}\xspace}
\newcommand{\TUCC}{\textsc{tucc}\xspace}
\newcommand{\cTUCC}{\textsc{closed tucc}\xspace}

\newcommand{\gTCC}{\textsc{(closed) tcc}\xspace}
\newcommand{\gTUCC}{\textsc{(closed) tucc}\xspace}

\begin{document}
\maketitle
\begin{abstract}
    A temporal (directed) graph is a graph whose edges are available only at specific times during its lifetime, $\tau$. Paths are sequence of adjacent edges whose appearing times are either strictly increasing or non-strictly increasingly (i.e., non-decreasing) depending on the scenario. Then, the classical concept of connected components and also of unilateral connected components in static graphs naturally extends to temporal graphs. 
    In this paper, we answer the following fundamental questions in temporal graphs. (i) What is the complexity of deciding the existence of a component of size $k$, parameterized by $\tau$, by $k$, and by $k+\tau$? We show that this question has a different answer depending on the considered definition of component and whether the temporal graph is directed or undirected.
    (ii) What is the minimum running time required to check whether a subset of vertices are pairwise reachable? A quadratic algorithm is known but, contrary to the static case, we show that a better running time is unlikely unless \textsc{SETH} fails. (iii) Is it possible to verify whether a subset of vertices is a component in polynomial time? We show that depending on the definition of temporal component this test is \NP-complete.
\end{abstract}

\section{Introduction}

A \emph{(directed) temporal graph $(G,\lambda)$ with lifetime $\tau$} consists of a (directed) graph $G$ together with a \emph{time-function} $\lambda: E(G)\to 2^{[\tau]}$ which tells when each edge $e\in E(G)$ is available along the discrete time interval $[\tau]$. Given $i\in [\tau]$, the \emph{snapshot $G_i$} refers to the subgraph of $G$ containing exactly the edges available in time $i$. Temporal graphs, also appearing in the literature under different names~\cite{casteigts2012time,BorgnatFGMRS07,LVM.18}, have attracted a lot of attention in the past decade, as many works have extended classical notions of Graph Theory to temporal graphs (we refer the reader to the surveys~\cite{LVM.18,M15} and the seminal paper~\cite{KKK00}).

A crucial characteristic of temporal graphs is that a $u,v$-walk/path in $G$ is valid only if it traverses a sequence of adjacent edges $e_1, \ldots, e_k$ at non-decreasing times $t_1\le \ldots\le t_k$, respectively, with $t_i\in \lambda(e_i)$ for every $i\in [k]$. Similarly, one can consider strictly increasing sequences, i.e. with $t_1<\ldots<t_k$. The former model is referred to as \emph{non-strict} model, while the latter as \emph{strict}. In both settings, we call such sequence a \emph{temporal $u,v$-walk/path}, and we say that \emph{$u$ reaches $v$}. For instance, in Figure~\ref{fig:example}, both blue and green paths are valid in the non-strict model, but only the green one is valid in the strict model, as the blue one traverses two edges with label $2$. The red path is not valid in both models.

\begin{figure}[t]
  \begin{center}
  \scriptsize
  \begin{subfigure}[b]{0.43\textwidth}
  \begin{tikzpicture}
 \vertex[label=left:$a$] (a) at (0,0) {};
 \vertex[label=left:$b$]  (b) at (0,2) {};
 \vertex[label=above:$c$]  (c) at (2,2) {};
 \vertex[label=below:$d$]  (d) at (2,0) {};
 \vertex[label=right:$e$]  (e) at (4,2) {};
 \vertex[label=right:$f$] (f) at (4,0) {};
 \vertex[label=below:$g$] (g) at (3,.75){};

 \draw[blue,thick,->] (a) to[black] node [left] {$\textcolor{blue}{1},5$} (b) ;
 \draw[->] (b)to[out=-60,in=60] node[right] {$1$}(a);
 \draw[blue,thick,->] (b)to[black]  node[anchor=south]{$\textcolor{blue}{2},6$}(c);
 \draw[blue,thick,->] (c)to node[black,above] {$\textcolor{blue}{2}$} (e);
 \draw[->] (e)to[out=210, in=-30 ] node[below] {$2$}(c);
 \draw[red,thick,->] (c) to[black] node[left]{$1,\textcolor{red}{3}$} (d);
 \draw[red,thick,->] (d) to[black] node[below] {$\textcolor{red}{1},2$} (f);
 \draw[green,thick,->] (f) to[black] node[sloped,above] {$\textcolor{green}{2}$} (g);
\draw[->] (f) to[black] node[right] {$3$} (e);
\draw[green,thick,->] (g) to[black] node[sloped, above]{$\textcolor{green}{3}$} (d);
 \draw[green,thick,->] (d) to[black] node[below]{$\textcolor{green}{4}$} (a);
  \end{tikzpicture}
  \vspace{-1.5cm}
  \end{subfigure}~\begin{subfigure}[b]{0.57\textwidth}
  \begin{tabular}{p{6.5cm}}
    $A'=\{a,b\}$ is a closed connected set, as $a$ and $b$ reach each other without using external vertices.\\
    \hline
    $A=\{a,b,c,d\}$ is a maximal closed connected set, i.e. a \cTCC.\\
    \hline
    $B=\{a,b,c,d,e\}$ is a \cTUCC but not a \cTCC as, using only vertices in $B$, $a,b,c,d$ reach each other, $e$ reaches all the vertices in $B$ and vice versa, except for $d$, which does not reach $e$. $B$ is also a \TCC, as $d$ can reach $e$ using the external vertex $f$.\\
    \hline
    $C=\{a,b,c,d,e,f\}$ is a \TUCC as $B$ forms a \cTUCC, $f$ is able to reach every other vertex directly or via the external vertex $g$. However, $C$ is not a \TCC as $a,b,c,e$ cannot reach $f$.
    \end{tabular}
    \end{subfigure}
  \end{center}
  \caption{On the left a temporal graph, where on each edge $e$ we depict $\lambda(e)$. Some of its components according to the non-strict model are reported on the right.
  \vspace{-0.5cm}}
\label{fig:example}
\end{figure}

The non-strict model is more appropriate in 
situations where the time granularity is relatively big.
This is the case in a disease-spreading scenario~\cite{zschoche2020complexity}, where the spreading speed might be unclear or in ``time-varying graphs'', as in in~\cite{nicosia.etal.12}, where a single snapshot corresponds to all the edges available in a time interval, e.g. the set of all the streets available in a day.
As for the strict model, it can represent the connections of the public transportation network of a city which are available only at precise scheduled times. All in all, there is a rich literature on both models (see \cite{abs-2202-12055,casteigts2021finding,rymar2021towards,zschoche2020complexity,haag2022feedback}), 
and this is why we explore both settings.

\textbf{Connected Sets and Components.} 
Given a temporal graph ${\cal G} = (G,\lambda)$, we say that $X\subseteq V(G)$ is a \emph{temporal connected set} if $u$ reaches $v$ \emph{and} $v$ reaches $u$, for every $u,v\in X$. 
Extending the classical notion of connected components in static graphs, in~\cite{BF.03} the authors define a \emph{temporal connected component} ({\TCC} for short) as a maximal connected set of ${\cal G}$. Such constraint can be strengthened to 
the existence of such paths using only vertices of $X$. Formally, $X$ is 
a \emph{closed temporal connected component} ({\cTCC} for short) if, for every $u,v\in X$, we have that $u$ reaches $v$ \emph{and} $v$ reaches $u$ through temporal paths contained in $X$. See Figure~\ref{fig:example} for an example of {\TCC} and {\cTCC}.

\textbf{Unilateral Connected Components.}
In the same fashion, also the concept of \emph{unilateral connected components} can be extended to temporal graphs. In static graph theory, they are a well-studied relaxation of connected components which asks for a path from $u$ to $v$ \underline{or} vice versa, for every pair $u,v$ in the component~\cite{A.76,B.72}. More formally, in a directed graph $G$, we say that $X\subseteq V(G)$ is a  \emph{unilateral connected connected set} if either $u$ reaches $v$ \underline{or} $v$ reaches $u$, for every $u,v\in X$. $X$ is a \emph{unilateral connected component} if it is maximal. In this paper, we introduce the definition of a \emph{(closed) unilateral temporal connected set/component}, which can be seen as the immediate translation of unilateral connected component to the temporal context. Formally, $X\subseteq V(G)$ is a \emph{temporal unilateral connected set} if $u$ reaches $v$ \emph{or} $v$ reaches $u$, for every $u,v\in X$, and it is a \emph{closed unilateral connected set} if this holds using paths contained in $X$. Finally, a \emph{(closed) temporal unilateral connected component} (\gTUCC for short) is a maximal (closed) temporal unilateral connected set. 
See again Figure~\ref{fig:example} for an example.

\paragraph{\textbf{Problems.}}
In this paper, we deal with four different definitions of temporal connected components, depending on whether they are unilateral or not, and whether they are closed or not. 
In what follows, we pose three questions, and we comment on partial knowledge about each of them. Later on, we discuss our results, which close almost all the gaps found in the literature. 
We start by asking the following. 

\begin{question}[Parameterized complexity] \label{quest:one}Deciding the existence of temporal components of size at least $k$ parameterized by  \emph{(i)} $\tau$, i.e. the lifetime,  \emph{(ii)} $k$, and \emph{(iii)} $k+\tau$.
\end{question}

In order to answer Question~1 for the strict model, there is a very simple parameterized reduction from $k$-clique, known to be \W[1]-hard when parameterized by $k$~\cite{downey1995fixed}, to deciding the existence of connected components (both closed or not and both unilateral or not) of size at least $k$ in undirected temporal graphs. This reduction has appeared in~\cite{C.18}. 
Given an undirected graph $G$, we can simply consider the temporal graph ${\cal G} = (G,\lambda)$ where $\lambda(uv)=\{1\}$ for all $uv\in E(G)$ (i.e., ${\cal G}$ is equal to $G$ itself). As $u$ temporally reaches $v$ if and only if $uv\in E(G)$, one can see that all those problems are now equivalent to deciding the existence of a $k$-clique in $G$. Observe that we get $\W[1]$-hardness when parameterized by $k$ or $k+\tau$, and para-$\NP$-completeness when parameterized by $\tau$, both in the undirected and the directed case.\footnote{In the directed case, it suffices to replace each edge of the input graph with two opposite directed edges between the same endpoints.} However, this reduction does not work in the case of the \emph{non-strict} model, leaving Question~\ref{quest:one} open. Indeed the reductions in~\cite{BF.03} and in~\cite{CCS.22} for {\gTCC}s, which work indistinctly for both the strict or the non-strict models, are not parameterized reductions. We also observe that the aforementioned reductions work on the non-strict model only for $\tau\ge 4$.

Another question of interest is the following. Letting $n$ be the number of vertices in ${\cal G}$ and $M$ be the number of \emph{temporal edges},\footnote{$M=\sum_{e\in E({\cal G})}|\lambda(e)|$.} it is known that, in order to verify whether $X\subseteq V(G)$ is a connected set in ${\cal G}$, we can simply apply $O(n)$ single source ``best'' path computations (see e.g.~\cite{CalamaiCM22,wu2014path}), resulting in a time complexity of $O(n\cdot M)$. This is $O(M^2)$ if ${\cal G}$ has no isolated vertices, a natural assumption when dealing with connectivity problems. As in static graphs testing connectivity can be done in linear time~\cite{hopcroft1973algorithm}, we ask whether the described algorithm can be improved.

\begin{question}[Lower bound on checking connectivity] \label{quest:lower}
Given a temporal graph ${\cal G}$ and a subset $X\subseteq V({\cal G})$, what is the minimum running time required to check whether $X$ is a (unilateral) connected set?
\end{question}

Finally we focus on one last question. 

\begin{question}[Checking maximality] \label{quest:maximal} Given a temporal graph ${\cal G}$ and a subset $X\subseteq V({\cal G})$, is it possible to verify, in polynomial time, whether $X$ is a component, i.e. a maximal (closed) (unilateral) connected set?\end{question}

For Question~\ref{quest:maximal}, we first observe that the property of being a temporal (unilateral) connected set is hereditary (forming an independence system~\cite{lawler1980generating}, see~\cite{ConteGMV19} for a survey about set systems), meaning that every subset of a (unilateral) connected set is still a (unilateral) connected set. For instance, in Figure~\ref{fig:example}, every subset of the connected set $B=\{a,b,c,d,e\}$ is a connected set. 
Also, checking whether $X'\subseteq V(G)$ is a temporal (unilateral) connected set can be done in time $O(n\cdot M)$, as discussed above. We can then check whether $X$ is a maximal such set in time $O(n^2\cdot M)$: it suffices to test, for every $v\in V(G)\setminus X$, whether by adding $v$ to $X$ we still get a temporal (unilateral) connected set. 
On the other hand, \emph{closed} connected (unilateral) sets are not hereditary, because by removing vertices from the set we could destroy the paths between other members of the set. This is the case for the closed connected set $A=\{a,b,c,d\}$ in Figure~\ref{fig:example}, since by removing $d$ there are no temporal paths from $c$ to $a$ nor $b$ anymore. 
This implies that the same approach as before does not work, i.e., we cannot check whether $X$ is maximal by adding to $X$ a \emph{single} vertex at a time, then checking for connectivity. 
For instance, the closed connected set $A'=\{a,b\}$ in Figure~\ref{fig:example} cannot be grown into the closed connected set $A$ by adding one vertex at a time, since both $A'\cup \{c\}$ and $A'\cup\{d\}$ are not closed connected sets. 
Hence, the answer to Question~\ref{quest:maximal} for closed sets does not seem easy, and until now was still open. 

We remark the important practical consequences of the latter question. Indeed, in practice, when trying to find structures of maximum size, a usual viable strategy is modifying \emph{backtracking listing algorithms} for efficient generating maximal structures (eventually with pruning strategies) and choosing the largest structures found~\cite{eblen2012maximum,bron1973algorithm}. Such algorithms typically solve the so-called \emph{extension problem}, that is generating all (or some of) the maximal solutions enlarging a partial one~\cite{kante2015polynomial,ConteGKMUW19,avis1996reverse,ConteGMV19}. Question~\ref{quest:maximal} implicitly asks whether efficient generation of {\cTCC}s or {\cTUCC}s is likely to exist or not.

\begin{table}[h]
    \centering\scriptsize
    \begin{tabular}{|c|c|c|c|}

    \hline
            & \textsc{Par. $\tau$} & \textsc{Par. $k$} & \textsc{Par. $k+\tau$} \\
    \hline
    \hline
    \multirow{2}{*}{\TCC}  &  \multirow{8}{*}{p-$\NP$ $\tau\ge 2$ (Th.~\ref{thm:param_tau})}  & $\W[1]$-h Dir. $\tau\geq 2$ (Th.~\ref{thm:W1h_dir_tau2}) & \\
    & & and Undir. (Th.~\ref{thm:TCCbyk}) & $\W[1]$-h Dir. (Th.~\ref{thm:W1h_dir_tau2})\\
    \cline{1-1}\cline{3-3}
    \multirow{2}{*}{\TUCC}  &  & $\W[1]$-h Dir. $\tau\geq 2$ (Th.~\ref{thm:W1h_dir_tau2}) & $\FPT$ Undir. (Th.~\ref{thm:FPT_TCC_undir})\\
              &  &  $\FPT$ Undir. (Th.~\ref{thm:FPT_TCC_undir}) & \\
    \cline{1-1}\cline{3-4}
    \multirow{2}{*}{\cTCC}  & & \multirow{2}{*}{$\W[1]$-h Dir. $\tau\ge 3$ (Th.~\ref{thm:cTCC_cTUCC_Whard})}  & \\
    &  & & $\W[1]$-h Dir. (Th.~\ref{thm:cTCC_cTUCC_Whard})\\
    \cline{1-1}\cline{3-3}
    %
    \multirow{2}{*}{\cTUCC}  &  & $\W[1]$-h Dir. $\tau\ge 3$ (Th.~\ref{thm:cTCC_cTUCC_Whard}) & $\FPT$ Undir. (Th.~\ref{thm:FPT_TCC_undir})\\
    &  & $\FPT$ Undir. (Th.~\ref{thm:FPT_TCC_undir})  &  \\
    \hline
    \end{tabular}
    \caption{A summary of our results for the parameterized complexity of computing components of size at least $k$ of a temporal graph ${\cal G}$ having lifetime $\tau$ in the \emph{non-strict} model. ``$\W[1]$-h'' stands for $\W[1]$-hardness and ``p$\NP$'' stands for para-\NP-completeness. For the strict model the entries are $\W[1]$-h in the third and fourth columns and p-\NP\ in the second one already for $\tau=1$, both for the directed and the undirected case.\vspace{-1cm}}
    \label{tab:results}
    
\end{table}

\paragraph{\textbf{Our results.}}
Our results concerning Question 1 are reported in Table~\ref{tab:results} for the non-strict model, since for the strict model all the entries would be \W[1]-hard or para-\NP-complete already for $\tau=1$, as we argued before. In the non-strict model, we observe instead that the situation is much more granulated. If $\tau=1$, then all the problems become the corresponding traditional ones in static graphs, which are all polynomial (see Paragraph ``Related works'').  
As for bigger values of $\tau$, the complexity depends on the definition of component being considered, and whether the temporal graph is directed or not. Table~\ref{tab:results} considers $\tau>1$, reporting on negative results ``$\tau\geq x$" for some $x$ meaning that the negative result starts to hold for temporal graphs of lifetime at least $x$.

The second column of Table~\ref{tab:results} addresses Question~\ref{quest:one}\emph{(i)}, i.e., parameterization by $\tau$. We prove that, for all the definitions of components being considered, the related problem becomes immediately para-\NP-complete as soon as $\tau$ increases from~1 to~2; this is done in Theorem~\ref{thm:param_tau}. This reduction improves upon the reduction of~\cite{BF.03}, which holds only for $\tau\geq 4$.

Question~\ref{quest:one}\emph{(ii)} (parameterization by $k$) is addressed in the third column of Table~\ref{tab:results}. Considering first directed temporal graphs, we prove that all the problems are $\W[1]$-hard. In particular, deciding the existence of a {\TCC} or \TUCC of size at least $k$ is $\W[1]$-hard already for $\tau\geq 2$ (Theorem~\ref{thm:W1h_dir_tau2}). As for the existence of closed components, $\W[1]$-hardness also holds as long as $\tau\ge 3$  (Theorem~\ref{thm:cTCC_cTUCC_Whard}). Observe that, since $\tau$ is constant in both results, these also imply the $\W[1]$-hardness results presented in the last column, thus answering also Question~\ref{quest:one}\emph{(iii)} (parameterization by $k+\tau$) for directed graphs. 
On the other hand, if the temporal graph is undirected, then the situation is even more granulated. 
Deciding the existence of a {\TCC} of size at least $k$ remains $\W[1]$-hard, but only if $\tau$ is unbounded. 
This is complemented by the answer to Question~\ref{quest:one}\emph{(iii)}, presented in the last column of Table~\ref{tab:results}: {\TCC} and (even) {\cTCC} are $\FPT$ on undirected graphs when parameterized by $k+\tau$ (Theorem~\ref{thm:FPT_TCC_undir}). 
We also give $\FPT$ algorithms when parameterized by $k$ for unilateral components, namely {\TUCC} and {\cTUCC}. Observe how this differs from {\TCC}, whose corresponding problem 
is $\W[1]$-hard, meaning that unilateral and traditional components behave very differently when parameterized by $k$.

In summary, Table~\ref{tab:results} answers Question~\ref{quest:one} for almost all the definitions of components, both for directed and undirected temporal graphs, leaving open only the problems of, given an undirected temporal graph, deciding the existence of a {\cTCC} of size at least $k$ when parameterized by $k$, and solving the same problem for {\cTCC} and {\cTUCC} in directed temporal graphs where $\tau=2$.

Concerning Questions~\ref{quest:lower} and~\ref{quest:maximal}, our results are summarized in Table~\ref{tab:res2}. All these results hold both for the strict and the non-strict models. For Question~\ref{quest:lower}, we prove that the trivial $O(M^2)$ algorithm to test whether $S$ is a (closed) (unilateral) connected set is best possible, unless the Strong Exponential Time Hypothesis (\textsc{SETH}) fails~\cite{impagliazzo2001complexity}. For Question~\ref{quest:maximal}, in the case of \TCC and \TUCC, we have already seen that checking whether a set $X\subseteq V$ is a component can be done in $O(n^2\cdot M)$. Interestingly, for \cTCC and \cTUCC, we answer negatively (unless \P=\NP) to Question~\ref{quest:maximal}.

\begin{table}[h]
    \centering\scriptsize
    \begin{tabular}{|c|c|c|c|}
 \hline
 & \textsc{Check whether $X\subseteq V$ is} & \textsc{Check whether $X\subseteq V$ is} \\
 & \textsc{a connected set} &  \textsc{a component} \\
 \hline
 \hline
{\TCC} & \multirow{4}{*}{$\Theta(M^2)$~(Th.~\ref{thm:lowerbound_algorithm})}  & \multirow{2}{*}{$O(n^2\cdot M)$} \\
\cline{1-1}
{\TUCC} & & \\
\cline{1-1}\cline{3-3}
{\cTCC} & & \multirow{2}{*}{\NP-c (Th.~\ref{thm:maximality})}\\
\cline{1-1}
{\cTUCC} & & \\
\hline
    \end{tabular}
\caption{Our results for Question~\ref{quest:lower} and Question~\ref{quest:maximal}, holding for both the \emph{strict} and the \emph{non-strict} models. Recall that a component is a (inclusion-wise) maximal connected set. The $O(\cdot)$ result is easy and explained in the introduction. $M$ (resp. $n$) denotes the number of temporal edges (resp. nodes) in ${\cal G}$.\vspace{-1cm}}
\label{tab:res2}
\end{table}

\medskip


\paragraph{\textbf{Related work}.}\label{sec:related}

The known reductions for temporal connected components in the literature~\cite{CCS.22,BF.03} which considers the non-strict setting are not parameterized and leave open the case when $\tau=2$ or $3$. The reductions we give here are parameterized (Theorems~\ref{thm:W1h_dir_tau2} and~\ref{thm:cTCC_cTUCC_Whard}) and Theorem~\ref{thm:param_tau_TUCC} closes also the cases $\tau=2$ and $3$. Furthermore, in~\cite{CCS.22} they show a series of interesting transformations but none of them allows us to apply known negative results for the strict model to the non-strict one.

There are many other papers about temporal connected components in the literature, including~in \cite{C.18}, where they give an example where there can be an exponential number of temporal connected components in the strict model.
In~\cite{nicosia.etal.12}, the authors show that the problem of computing {\TCC}s is a particular case of finding cliques in the so-called \emph{affine graph}. This does not imply that the problem is $\NP$-complete as claimed.
Further related works include recent papers on giant components and connectivity in randomized temporal graphs~\cite{abs-2205-14888,CasteigtsRRZ21} and on networks with continuous varying-time~\cite{AS.19}.

Other notions of temporal components in the literature include \emph{temporal out-component} (resp. \emph{in-component}) in~\cite{NTMMRL.13}, $\Delta$-component in~\cite{GCLL.15}, weakly connected components in~\cite{NTMMRL.13}. The latter applies if the temporal graph $(G,\lambda)$ is directed and, as in the analogue case in static directed graphs, it simply ignore directions of the edges and consider the undirected version of the underlying graph $G$. In this paper we implicitly consider also weakly temporal connected components as studying how to compute them in directed temporal graphs is the same as studying how to compute {\gTCC}s in undirected temporal graphs.

Finally, we remark that there are many results in the literature concerning unilateral components in static graphs, also with applications to community detection~\cite{LP.11}.
Even though the number of unilateral components in a graph is exponential~\cite{A.76}, deciding whether there is one of size at least $k$ is polynomial. In~\cite[Theorem 3]{A.76}, they prove that this corresponds to deciding whether a DAG with weights on its vertices has a path of weight at least $k$, which in turn can be done in polynomial time by slightly modifying the algorithm for longest paths~\cite{book_Sedgewick_Wayne}. Additionally, in~\cite{A.76,Cheston78}, the authors propose a listing algorithm, and in~\cite{BG.08}, a characterization of unilaterally connected graphs is presented. 
Further related works include~\cite{MS.10}, and ~\cite{FL.80}.
There were no results about unilateral components in temporal graphs until now.

\paragraph{\textbf{Preliminaries}.} For further formal definitions, we refer to Appendix~\ref{app:prelim}.  For basic graph theory concepts and notation, we refer to~\cite{West.book}, and to~\cite{DF13,CyganFKLMPPS15} for basic background on parameterized complexity. 

\paragraph{\textbf{Structure of the Paper}.} In what follows, we  present sketches of proofs of all results, with the results concerning Question~\ref{quest:one} presented in Section~\ref{sec:param_results}, and the ones related to Questions~\ref{quest:lower} and~\ref{quest:maximal} presented in Section~\ref{sec:lowerbound_maximality}. In Section~\ref{sec:conclusion}, we present our concluding remarks.

\section{Parameterized Complexity Results}\label{sec:param_results}

This section is devoted to answer Question~\ref{quest:one} and prove the results summarized in Table~\ref{tab:results}.

\paragraph{Parameterization by $\tau$.} We start by proving the result in the first column of Table~\ref{tab:results} about para-\NP-completeness wrt the lifetime $\tau$, which applies to all the definitions of components. For \gTCC, we do a reduction from  \textsc{Maximum Edge Biclique}  (\textsc{MEBP} for short), which consists in, given a bipartite graph $G$ and an integer $k$, deciding whether $G$ has a biclique with at least $k$ edges. It was proved to be {\NP}-complete in~\cite{P.03}. 
Using the same construction, we prove hardness of \gTUCC reducing from \textsc{$2K_2$-free Edge Subgraph}, which consists in, given a bipartite graph $G$ and an integer $k$, deciding whether $G$ has a $2K_2$-free subgraph with at least $k$ edges.
This was proved to be $\NP$-complete in~\cite{Y.81}. 

The main idea of the reductions is to generate a temporal graph $\mathcal{G}$ whose base graph is the line graph $L$ of a bipartite graph $H$ with parts $X,Y$.
We make active in timestep $1$ every edge of a clique in $L$ related to vertices in $X$, and in timestep $2$ every edge of a clique related to vertices in $Y$.
Doing so, we ensure that any pair of vertices of $\mathcal{G}$ associated with a biclique in $H$ reach one another in $\mathcal{G}$. We prove that there exists a biclique in $H$ with at least $k$ edges if and only if there exists a {\cTCC} in $\mathcal{G}$ of size at least $k$. The result extends to {\TCC}s, as every {\TCC} is also a {\cTCC}. For the unilateral case, we can relax the biclique to a $2K_2$-free graph since only one reachability relation is needed.
As a result, we get the following, whose formal constructions and correctness are proven in Appendix~\ref{app:param_tau}.

\begin{theorem}
\label{thm:param_tau}
For every fixed $\tau \geq 2$ and given a temporal graph $\mathcal{G}$ of lifetime $\tau$ and an integer $k$, it is \NP-complete to decide if $\mathcal{G}$ has a {\gTCC} or a {\gTUCC} of size at least $k$, even if the base graph of $\mathcal{G}$ is the line graph of a bipartite graph.
\label{thm:param_tau_TUCC}
\end{theorem}

\paragraph{\W[1]-hardness by $k$.} We now focus on proving the \W[1]-hardness results in the second column of Table~\ref{tab:results} concerning parameterization by $k$, which also imply some of the results of the third column. The following \W[1]-hardness results (Theorem~\ref{thm:TCCbyk},~\ref{thm:cTCC_cTUCC_Whard}, and~\ref{thm:W1h_dir_tau2}) are parameterized reductions from $k$-\textsc{Clique}. 
The general objective is constructing a temporal graph $\mathcal{G}$ in a way that vertices in  $\mathcal{G}$ are in the same component if and only if the corresponding nodes in the original graph are adjacent. 
Notice that we have to do this while: (i) ensuring that the size of the desired component is $f(k)$ for some computable function $k$ (i.e., this is a parameterized reduction); and (ii) avoiding that the closed neighborhood of a vertex forms a component, so as to not a have a false ``yes'' answer to $k$\textsc{Clique}. 
To address these tasks, we rely on different techniques. 
The first reduction concerns \TCC in undirected graphs and requires $\tau$ to be unbounded, as for $\tau$ bounded we show that the problem is {\FPT} by $k+\tau$ (Theorem~\ref{thm:FPT_TCC_undir}). The technique used is a parameterized evolution of the so-called \emph{semaphore} technique used in~\cite{BF.03,CCS.22}, which in general replaces edges by labeled diamonds to control paths of the original graph. 
However, while the original reduction gives labels in order to ensure that paths longer than one are broken, the following one allows the existence of paths longer than one. But if a temporal path from $u$ to $v$ exists for $uv\notin E(G)$, then the construction ensures the non-existence of temporal paths from $v$ to $u$. Because of this property, the reduction does not extend to {\TUCC}s, which we prove to be {\FPT} when parameterized by $k$ instead (Theorem~\ref{thm:FPT_TCC_undir}). The interested reader can skip directly to Appendix~\ref{app:TCC_Undirected_Whard} for the complete proof.

\begin{theorem}\label{thm:TCCbyk}
Given a temporal graph $\mathcal{G}$ and an integer $k$, deciding if $\mathcal{G}$ has a \TCC of size at least $k$ is $\W[1]$-hard with parameter $k$.
\end{theorem}

\begin{figure}[t]
\begin{center}
\scalebox{.6}{
\begin{tikzpicture}[scale=1,vertex/.style={circle, minimum size=0.2cm, draw, inner sep=1pt}, blackvertex/.style={draw,circle,minimum size=5pt,inner sep=0pt, fill=black}]

\begin{scope}[xshift = -3cm]
\draw[rounded corners] (0,0) rectangle  (1.2,5) node [above, xshift = -.5cm] {\large$V$};

\draw[rounded corners] (2,0) rectangle  (3.2,5) node [above, xshift = -.5cm] {\large$V'$};

\foreach \x/\i in {0.5/n, 2.5/3, 3.5/2, 4.5/1}{
\node[blackvertex, label = left:{\large$u_{\i}$}] (u\i) at (0.75, \x) {};
\node[blackvertex, label = right:{\large$u'_{\i}$}] (v\i) at (2.45, \x) {};
\draw[thick] (u\i) -- (v\i) node[midway, fill=white] {$0$};
}

\path (u3) -- (un) node [font=\Huge, midway, sloped] {$\dots$};
\path (v3) -- (vn) node [font=\Huge, midway, sloped] {$\dots$};
\end{scope}

\begin{scope}[xshift=6cm, yshift = 5.5cm]
\node[blackvertex, scale = 1.1, label={\large$u$}] (u) at (0,0) {};
\node[blackvertex, scale = 1.1, label={\large$v$}] (v) at (2,0) {};
\draw[thick] (u) -- (v) node[midway, above]{\large$e_i$};
\end{scope}

\begin{scope}[xshift = 6cm, yshift = 3.5cm]
\node[blackvertex, scale = 1.1, label=left:{\large$u$}] (u) at (-2,0) {};
\node[blackvertex, scale = 1.1, label=right:{\large$v$}] (v) at (4,0) {};
\node[blackvertex, scale = 1.1, label=below:{\large$h_{uv}$}] (huv) at (1,1) {};
\node[blackvertex, scale = 1.1, label={\large$h_{vu}$}] (hvu) at (1,-1) {};

\draw[thick] (u) -- (huv) node[midway, fill=white] {$i$};
\draw[thick] (hvu) -- (v) node[midway, fill=white] {$i$};
\draw[thick] (huv) -- (v) node[midway, fill=white] {$m+i$};
\draw[thick] (u) -- (hvu) node[midway, fill=white] {$m+i$};

\begin{scope}[yshift = -2.5cm]
\node[blackvertex, scale = 1.1, label=left:{\large$u'$}] (u1) at (-2,0) {};
\node[blackvertex, scale = 1.1, label=right:{\large$v'$}] (v1) at (4,0) {};
\node[blackvertex, scale = 1.1, label=below:{\large$h'_{uv}$}] (huv1) at (1,1) {};
\node[blackvertex, scale = 1.1, label={\large$h'_{vu}$}] (hvu1) at (1,-1) {};

\draw[thick] (u1) -- (huv1) node[midway, fill=white] {$2m+i$};
\draw[thick] (hvu1) -- (v1) node[midway, fill=white] {$2m+i$};
\draw[thick] (huv1) -- (v1) node[midway, fill=white] {$3m+i$};
\draw[thick] (u1) -- (hvu1) node[midway, fill=white] {$3m+i$};
\end{scope}
\draw[thick] (u) -- (u1) node [midway, fill=white] {$0$};
\draw[thick] (v) -- (v1) node [midway, fill=white] {$0$};

\end{scope}

\end{tikzpicture}
}
\end{center}
\caption{Construction used in the proof of Theorem~\ref{thm:TCCbyk}. On the left, the two copies of $V(G)$ and the edges between them, active in timestep $0$. On the right, the edge $e_i \in E(G)$ and the associated gadget in $\mathcal{G}$.}
\label{fig:theorem-2}
\end{figure}
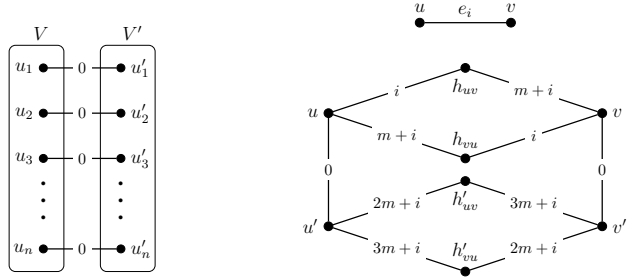

\begin{proof}
We make a parameterized reduction from $k$-\textsc{Clique}. 
Let $G$ be graph and $k \geq 3$ be an integer.
We construct the temporal graph $\mathcal{G} = (G', \lambda)$ as follows. See Figure~\ref{fig:theorem-2} to follow the construction. 
First, add to $G'$ every vertex in $V(G)$ and make $V = V(G)$.
Second, add to $G'$ a copy $u'$ of every vertex $u \in V$ and define $V' = \{u' \mid u \in V\}$.
Third, for every pair $u,u'$ with $u\in V$ and $u' \in V'$ add the edge $uu'$ to $G'$ and make all such edges active at timestep $0$.
Fourth, consider an arbitrary ordering $e_1, \ldots, e_m$ of the edges of $G$ and,
for each edge $e_i = uv$, create four new vertices $\{h_{uv},h_{vu},h'_{uv},h'_{vu}\mid uv\in E(G)\}$, adding edges: 
 \begin{itemize}
  \item $uh_{uv}$ and $vh_{vu}$, active at time $i$;
  \item $u'h'_{uv}$ and $v'h'_{vu}$, active at time $2m+i$;
  \item $h_{vu}u$ and $h_{uv}v$, active at time $m+i$; and
  \item $h'_{vu}u'$ and $h'_{uv}v'$, active at time $3m+i$.
\end{itemize}
Denote the set $\{h_{uv},h_{vu}\mid uv\in E(G)\}$ by $H$, and the set $\{h'_{uv},h'_{vu}\mid uv\in E(G)\}$  by $H'$. 
We now prove that $G$ has a clique of size at least $k$ if and only if ${\cal G}$ has a \TCC of size at least $2k$.
Given a clique $C$ in $G$, it is easy to check that $C\cup\{u\in V'\mid u\in C\}$ is a \TCC, and because of space constraints we present the formal argument only in Appendix~\ref{app:TCC_Undirected_Whard}. 

Now, let $S\subseteq V(G')$ be a \TCC of ${\cal G}$ of size at least $2k$. 
We want to show that either $C = \{u\in V(G)\mid u\in S\cap V\}$ or $C'= \{u\in V(G)\mid u'\in S\cap V'\}$ is a clique of $G$ of size at least $k$. 
This part of the proof combines a series of useful facts, which we cannot include here due to space constraints. The full proof can be found Appendix~\ref{app:TCC_Undirected_Whard}, and in what follows we present a sketch of it.

First, we argue that both $C$ and $C'$ are cliques in $G$.
Then, by observing that the only edges between $V \cup H$ and $V' \cup H'$ are those incident to $V$ and $V'$ at timestep $0$, we conclude that either $S \subseteq V\cup H$ or $S \subseteq V' \cup H'$.
Since the cases are similar, we assume the former.
If $|S\cap V| \geq k$, then $C$ contains a clique of size at least $k$ and the result follows.
Otherwise, we define $E_S = \{uv \in E(G) \mid \{h_{uv},h_{vu}\} \cap S \neq \emptyset\}$.
That is, $E_S$ is the set of edges of $G$ related to vertices in $S \cap H$.
We then prove the following claim.
\begin{clm*}\label{claim:sketch-thm2}
\emph{Let $a,b \in S \cap H$ be associated with distinct edges $g,g'$ of $G$ sharing an endpoint $v$.
If $u$ and $w$ are the other endpoints of $g$ and $g'$, respectively, then $u$ and $w$ are also adjacent in $G$.
Additionally, either $\lvert S\cap \{h_{xy},h_{yx}\}\rvert \le 1$ for every $xy\in E(G)$, or $\lvert S\cap H\rvert \le 2$.}
\end{clm*}

To finish the proof, we first recall that we are in the case $\lvert S\cap H\rvert \ge k+1$.
By our assumption that $k\ge 3$, note that the above claim gives us that $\lvert S\cap \{h_{xy},h_{yx}\}\rvert \le 1$ for every $xy\in E(G)$, which in turn implies that $\lvert E_S\rvert = \lvert S\cap H\rvert$. Additionally, observe that, since $\lvert S\cap H\rvert \ge 4$, the same claim also gives us that there must exist $w\in V$ such that $e$ is incident to $w$ for every $e\in E_S$. 
Indeed, the only way that $3$ distinct edges can be mutually adjacent without being all incident to a same vertex is if they form a triangle.
Supposing that $3$ edges in $E_S$ form a triangle $T = (a,b,c)$, since $\lvert E_S\rvert \ge 4$, there exists an edge $e\in E_S\setminus E(T)$.
But now, since $G$ is a simple graph, $e$ is incident to at most one between $a$, $b$ and $c$, say $a$.
We get a contradiction wrt the aforementioned claim as in this case $e$ is not incident to edge $bc\in E_S$. 
Finally, by letting $C'' = \{v_1,\ldots,v_k\}$ be any choice of $k$ distinct vertices such that $\{wv_1,\ldots,wv_k\}\subseteq E_S$, our claim gives us that $v_i$ and $v_j$ are adjacent in $G$, for every $i,j\in [k]$; i.e., $C''$ is a $k$-clique in $G$.
\end{proof}

The following result concerns \TCC and \TUCC in directed temporal graphs. It is important to remark that for \TCC and $\tau$ unbounded, we already know that the problem is \W[1]-hard because of Theorem~\ref{thm:TCCbyk} which holds for undirected graphs and extends to directed ones. However, the following reduction applies specifically for directed ones already for $\tau=2$. The technique used here is the previously mentioned semaphore technique, made parameterized by exploiting the direction of the edges. For the sake of space, the construction is shown only in Figure~\ref{fig:thm_3}, which shows how to obtain the temporal graph $\mathcal{G}$ in Figure~\ref{fig:thm_3}(b) from the graph $G$ in Figure~\ref{fig:thm_3}(a) in a way that, for every integer $k \geq 3$,  graph $G$ has a clique of size at least $k$ if and only if $\mathcal{G}$ has a \TCC of size at least $k$. The formal construction and proof are given in Appendix~\ref{app:W1h_dir_tau2}

\begin{theorem}\label{thm:W1h_dir_tau2}
Given a directed temporal graph $\mathcal{G}$ and an integer $k$, deciding if $\mathcal{G}$ has a \TCC of size at least $k$ is $\W[1]$-hard with parameter $k$, even if $\mathcal{G}$ has lifetime $2$. The same holds for \TUCC.
\end{theorem}
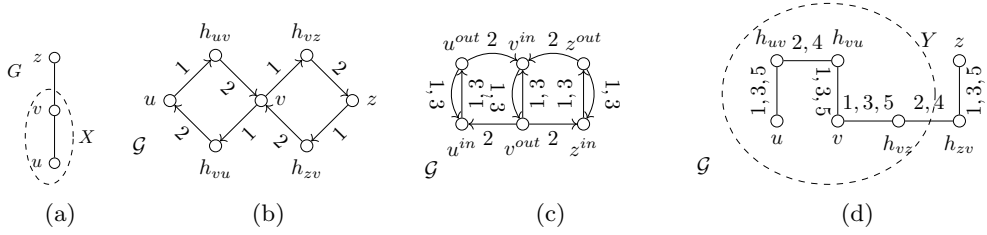
\begin{figure}[t]
    \centering
    \begin{subfigure}[b]{0.12\textwidth}
    \scalebox{0.7}{\begin{tikzpicture}[scale=1]

\node (G) at (-.75,1.75) {$G$};
\node[dashed,draw, ellipse, minimum width=.9cm, minimum height=1.8 cm] at (-.1,.5) {};
\node[label=right:$X$] at (0.2,.5) {};

\vertex[label=left:$u$] (u) at (0,0) {};
\vertex[label=left:$v$] (v) at (0,1) {};
\vertex[label=left:$z$] (z) at (0,2) {};

\draw (u) -- (v) -- (z);

\end{tikzpicture}}
    \caption{}
    \end{subfigure}~\begin{subfigure}[b]{0.29\textwidth}
    \scalebox{0.8}{\begin{tikzpicture}[scale=1]
 \node (G) at (-.5,-.75) {$\mathcal{G}$};
\vertex[label=left:$u$] (u) at (0,0) {};
\vertex[label=above:$h_{uv}$] (huv) at (.75,.75) {};
\vertex[label=below:$h_{vu}$] (hvu) at (.75,-.75) {};
\vertex[label=right:$v$] (v) at (1.5,0) {};
\vertex[label=above:$h_{vz}$] (hvz) at (2.25,.75) {};
\vertex[label=below:$h_{zv}$] (hzv) at (2.25,-.75) {};
\vertex[label=right:$z$] (z) at (3,0) {};

 \draw[->] (u) to node[scale=1,sloped,above] {$1$} (huv);
 \draw[->] (huv) to node[scale=1,sloped,below] {$2$}  (v);
 \draw[->] (v) to node[scale=1,sloped,below] {$1$} (hvu);
 \draw[->] (hvu) to node[scale=1,sloped,below] {$2$} (u);
 \draw[->] (v) to node[scale=1,sloped,above] {$1$} (hvz);
 \draw[->] (hvz) to node[scale=1,sloped,above] {$2$} (z);
 \draw[->] (z) to node[scale=1,sloped,below] {$1$} (hzv);
 \draw[->] (hzv) to node[scale=1,sloped,below,pos=.3] {$2$} (v);

\end{tikzpicture}}
    \caption{}
    \end{subfigure}~\begin{subfigure}[b]{0.27\textwidth}
    \scalebox{0.8}{\begin{tikzpicture}[scale=1]
 \node (G) at (-.5,-.75) {$\mathcal{G}$};
\vertex[label=below:$u^{in}$] (u-) at (0,0) {};
\vertex[label=above:$u^{out}$] (u+) at (0,1) {};
\vertex[label=above:$v^{in}$] (v-) at (1,1) {};
\vertex[label=below:$v^{out}$] (v+) at (1,0) {};
\vertex[label=below:$z^{in}$] (z-) at (2,0) {};
\vertex[label=above:$z^{out}$] (z+) at (2,1) {};

  \draw[->] (u-) to node[scale=1,sloped,anchor=north,below,pos=.5] {$1,3$} (u+);
  \draw[->] (u+) to[out=240,in=120] node[scale=1,sloped,anchor=east,below,pos=.5] {$1,3$}  (u-);
  \draw[->] (u+) to[out=30,in=150] node[scale=1,sloped,anchor=east,above,pos=.5] {$2$} (v-);
   \draw[->] (v+) to node[scale=1,sloped,anchor=west,below,pos=.55] {$\!2$} (u-);
  \draw[->] (v+) to node[scale=1,sloped,anchor=west,below,pos=.5] {$1,3$} (v-); 
  \draw[->] (v-) to[out=240,in=120] node[scale=1,sloped,below,pos=.5] {$~1,3$}  (v+);
   \draw[->] (z+) to[out=150,in=30] node[scale=1,sloped,anchor=east,above,pos=.5] {$2$} (v-);
  \draw[->] (v+) to node[scale=1,sloped,anchor=east,below,pos=.5] {$2$} (z-);
  \draw[->] (z-) to node[scale=1,sloped,anchor=west,above,pos=.5] {$1,3$} (z+);
  \draw[->] (z+) to[out=-60,in=60] node[scale=1,sloped,anchor=east,above,pos=.5] {$1,3$}  (z-);

\end{tikzpicture}}
    \caption{}
    \end{subfigure}~\begin{subfigure}[b]{0.34\textwidth}
    \scalebox{0.8}{\begin{tikzpicture}[scale=1]
 \node[dashed,draw,ellipse, minimum width=3.5cm,minimum height=3cm] at (.85,.45){};
 \node (G) at (-1.2,-.75) {$\mathcal{G}$};
\node (Y) at (2.5,1.35) {$Y$}{};
\vertex[label=below:$u$] (u) at (0,0) {};
\vertex[label=above:$\!\!\!h_{uv}$] (huv) at (0,1) {};
\vertex[label=above:$~~~h_{vu}$] (uvh) at (1,1) {};
\vertex[label=below:$v$] (v) at (1,0) {};
\vertex[label=below:$h_{vz}$] (hvz) at (2,0) {};
\vertex[label=below:$h_{zv}$] (vzh) at (3,0) {};
\vertex[label=above:$z$] (z) at (3,1) {};

\draw (u) to node[scale=1,anchor=south,sloped, above]{$1,3,5$}  (huv);
\draw (huv) to node[scale=1,sloped, above] {$2,4$} (uvh);
\draw (uvh) to node[scale=1,anchor=north,sloped, below]{$1,3,5$} (v);
\draw (v) to node[scale=1,sloped, above] {$1,3,5$} (hvz);
\draw (hvz) to node[scale=1,sloped, above] {$2,4$} (vzh);
\draw (vzh) to node[scale=1,anchor=north,sloped, below]{$1,3,5$} (z);

\end{tikzpicture}}
    \caption{}
    \end{subfigure}
    \caption{Examples for some of our reductions. Given the graph in (a), Theorem~\ref{thm:W1h_dir_tau2} constructs the directed temporal graph in (b), Theorem~\ref{thm:cTCC_cTUCC_Whard} constructs the directed temporal graph in (c), and, given additionally set $X$ in (a), Theorem~\ref{thm:maximality} contructs the temporal graph $\mathcal{G}$ and set $Y$ in (d).}
    \label{fig:thm_3}
\end{figure}

The next result concerns closed {\TCC}s and {\TUCC}s. In this case, we also reduce from $k$-\textsc{Clique}, but we cannot apply the semaphore technique as before. Indeed, as we are dealing with closed components, nodes must be reachable using vertices inside the components, while the semaphore technique would make them reachable via additional nodes, which do not necessarily reach each other. For this reason, in the following we introduce a new technique subdividing nodes, instead of edges, in order to break paths of the original graph of length longer than one, being careful to allow that these additional nodes reach each other. The construction is shown in Figure~\ref{fig:thm_3}, which shows how to construct temporal graph $\mathcal{G}$ in Figure~\ref{fig:thm_3}(c), given graph $G$ in Figure~\ref{fig:thm_3}(a) in a way that graph $G$ has a clique of size $k$ if and only if $\mathcal{G}$ has a \cTCC (\TUCC) of size at least $2k$. The formal construction and proof are given in Appendix~\ref{app:cTCC_cTUCC_Whard}. 

\begin{theorem}\label{thm:cTCC_cTUCC_Whard}
Given a directed temporal graph $\mathcal{G}$ and an integer $k$, deciding if $\mathcal{G}$ has a \cTCC of size at least $k$ is $\W[1]$-hard with parameter $k$, even if $\mathcal{G}$ has lifetime $3$. The same holds for \cTUCC.
\end{theorem}

\paragraph{FPT algorithms.} We now show our $\FPT$ algorithms to find {\gTCC}s and {\gTUCC}s in undirected temporal graphs, as for directed temporal graphs we have proved \W[1]-hardness. In particular, we prove the following result, whose proof is shown in Appendix~\ref{app:FPT_TCC_undir}.
\begin{theorem}\label{thm:FPT_TCC_undir}
Given a temporal graph ${\cal G} = (G,\lambda)$ on $n$ vertices and with lifetime $\tau$, and a positive integer $k$, there are algorithms running in time
\begin{enumerate}
    \item $O(k^{k\cdot \tau}\cdot n)$ that decides whether there is a \TCC of size at least $k$; 
    \item $O(2^{k^\tau}\cdot n)$ that decides whether there is a \cTCC of size at least $k$;
    \item $O(k^{k^2}\cdot n)$ that decides whether there is a \TUCC of size at least $k$; and
    \item $O(2^{k^k}\cdot n)$ that decides whether there is a \cTUCC of size at least $k$.
\end{enumerate}
\end{theorem}
\begin{proof}[Sketch] The \emph{reachability digraph} $R$ associated to $\mathcal{G}$ is a directed graph with the same vertex set as $\mathcal{G}$, and such that $uv$ is an edge in $R$ if and only $u$ reaches $v$ in $\mathcal{G}$ and $u\neq v$. This is related to the \emph{affine} graph in~\cite{nicosia.etal.12}. Observe that finding a \TCC (\TUCC) in $\mathcal{G}$ of size at least $k$ is equivalent to finding a set $S\subseteq V(\mathcal{G})$ in $R$ of size \emph{exactly} $k$ such that $uv\in E(R)$ and (or) $vu\in E(R)$ for every pair $u,v\in V(R)$. As for finding a \cTCC (\cTUCC), we need to have the same property, except that all subsets of size \emph{at least $k$} must be tested (recall that being a closed connected (unilateral) set is not hereditary). Therefore, if $\Delta = \Delta(R)$, then testing connectivity takes time $O(k^\Delta\cdot n)$ (it suffices to test all subsets of size $k-1$ in $N(u)$, for all $u\in V(R)$), while testing closed connectivity takes time $O(2^{\Delta}\cdot n)$ (it suffices to test all subsets of size \emph{at least} $k-1$ in $N(u)$, for all $u\in V(R)$). The proofs then consist in bounding the value $\Delta$ in each case. \qed
\end{proof}

It is important to observe that, for unilateral components, these bounds depend only on $k$, while for {\TCC}s and {\cTCC}s they depend on both $k$ and $\tau$. This is consistent with the fact that we have proved that for \TCC the problem is $\W[1]$-hard when parameterized just by $k$ (Theorem~\ref{thm:TCCbyk}).

\section{Checking Connectivity and Maximality}\label{sec:lowerbound_maximality}

This section is focused on Questions~\ref{quest:lower} and Question~\ref{quest:maximal}.
The former is open for all definitions of components for both the strict and the non-strict models. We answer to the question providing the following conditional lower bound, which holds for both models, where the notation $\Tilde{O}(\cdot)$ ignores polylog factors. Due to space constraints, its proof is shown in Appendix~\ref{app:lower}. We apply the technique used for instance in~\cite{BorassiCH16,puatracscu2010possibility,williams2010subcubic} to prove lower bounds for polynomial problems, falling within the fine-grained complexity framework.

\begin{theorem}\label{thm:lowerbound_algorithm}
Consider a temporal graph ${\cal G}$ on $M$ temporal edges. There is no algorithm running in time $\Tilde{O}(M^{2-\epsilon})$, for some $\epsilon$, that decides whether $G$ is  temporally (unilaterally) connected, unless \textsc{SETH} fails.
\end{theorem}

We now focus on Question~\ref{quest:maximal}. We prove the results in the second column of Table~\ref{tab:res2}, about the problem of deciding whether a subset of vertices $Y$ of a temporal graph is a component, i.e. a maximal connected set. The question is open both for the strict and the non-strict model. We argued already in the introduction that this is polynomial for \TCC and \TUCC for both models. In the following we prove $\NP$-completeness for \cTCC and \cTUCC on undirected graphs. The results extend to directed graphs as well.

\begin{theorem}
    Let ${\cal G}$ be a (directed) temporal graph, and $Y\subseteq V({\cal G})$. Deciding whether $Y$ is a \cTCC is $\NP$-complete. The same holds for \cTUCC.
    \label{thm:maximality}
\end{theorem}

\begin{proof}
We reduce from the problem of deciding whether a subset of vertices $X$ of a given a graph $G$ is a maximal 2-club, where a 2-club is a set of vertices $C$ such that $G[C]$ has diameter at most~2. This problem has been shown to be $\NP$-complete in~\cite{PajouhB12}. Let us first focus on the strict model. 
In this case, given $G$ we can build a temporal graph $\mathcal{G}$ with only two snapshots, both equal to $G$. Observe that $X$ is a 2-club in $G$ if and only if $X$ is a \cTCC in $\mathcal{G}$. Indeed, because we can take only one edge in each snapshot and $\tau=2$, we get that temporal paths will always have length at most 2. This also extends to \cTUCC by noting that all paths in $\mathcal{G}$ can be temporally traversed in both directions. 

In the case of the non-strict model, the situation is more complicated as in each snapshot we can take an arbitrary number of edges resulting in paths arbitrarily long. We show the construction for \cTCC in what follows, proving its correctness in Appendix~\ref{app:maximality}. 
Let $\mathcal{G}$ be obtained from $G$ by subdividing each edge $uv$ of $G$ twice, creating vertices $h_{uv}$ and $h_{vu}$, with $\lambda(uh_{uv}) = \lambda(vh_{vu}) = \{1,3,5\}$, and $\lambda(h_{uv}h_{vu}) = \{2,4\}$. 
See Figure~\ref{fig:thm_3} (d) for an illustration.

Given $(G,X)$, the instance of maximal $2$-club, we prove that $X$ is a maximal $2$-club in $G$ iff $Y = X\cup N_H(X)$ is a \cTCC in ${\cal G}$. For this, it suffices to prove that, given $X'\subseteq V(G)$ and defining $Y'$ similarly as before w.r.t. $X'$, we have that $G[X']$ has diameter at most~2 iff $Y'$ is a closed temporal connected set.
The proof extends to \cTUCC by proving that every \cTCC is also a \cTUCC and vice-versa.
\end{proof}

\section{Concluding remarks}\label{sec:conclusion}

In this paper, we revisit the notion of connected components in temporal graphs introduced in~\cite{BF.03} from the point of view of parameterized complexity.
We then consider unilateral connectivity in temporal graphs, and investigate all related problems, in both the strict and the non-strict setting, as well as both for directed and undirected temporal graphs, parameterizing by the size $k$ of the desired component, the lifetime $\tau$ of the considered (directed) temporal graph $\mathcal{G}$, and by $k+\tau$. We classify all possible entries in Table~\ref{tab:results}, leaving open just the following questions.
\begin{question}
   Given an undirected temporal graph $\mathcal{G}$, and considering parameterization by $k$, the size of the searched component, what is the complexity of deciding the existence of a \cTCC?
\end{question}

\begin{question}
   Given a directed temporal graph $\mathcal{G}$ with lifetime~2, and considering parameterization by $k$, the size of the searched component, what is the complexity of deciding the existence of a \textsc{Closed} \TCC (\TUCC)?
\end{question}

We additionally prove a lower bound for testing connectivity, and prove that deciding maximality of closed (unilateral) connectivity is $\NP$-complete.

\bibliographystyle{plain}
\bibliography{main}

\appendix

\section{Preliminaries and notation}
\label{app:prelim}

Given a graph $G = (V, E)$, directed or not, and a set $X \subseteq V(G)$ we write $G[X]$ for the subgraph of $G$ induced by $X$.
If $e$ is an edge of a directed or undirected graph with \emph{endpoints} $u$ and $v$, we may refer to $e$ as $(u,v)$ and say that $e$ is \emph{incident} to $u$ and $v$.
If $e$ is an edge\footnote{We refer to arcs of directed graphs as edges (following the notation in~\cite{West.book}} from $u$ to $v$ of a directed graph, we say that $e$ has \emph{tail} $u$, \emph{head} $v$, and is \emph{oriented} from $u$ to $v$.
The \emph{degree} $d_G(v)$ of a vertex $v$ of a (directed) graph $G$ is the number of edges of $G$ incident to $v$.
We denote by $\Delta(G)$ the maximum degree of a vertex of $G$.
The \emph{neighborhood} $N_G(v)$ of $v$ is the set $\{u \in V(G) \mid (u,v) \in E(G)\}$.
If $D$ is a digraph, he \emph{in-neighborhood} $N^-_D(v)$ of $v$ is the set $\{u \in V(D) \mid (u,v) \in E(D)\}$, and the \emph{out-neighborhood} $N^+_D(v)$ is the set $\{u \in V(D) \mid (v,u) \in E(D)\}$.

An undirected graph $G$ is said to be \emph{simple} if there is at most one edge between every pair of vertices of $G$ (i.e., there are no \emph{parallel} edges).
We say that $G$ is \emph{bipartite} if there is a partition of $V(G)$ into two non-empty sets $X$ and $Y$ such that every edge of $G$ has one endpoint in $X$ and the other endpoint in $Y$.
A \emph{matching} of $G$ is a set of edges $M \subseteq E(G)$ such that no two edges of $M$ share an endpoint; i.e., they all are pairwise \emph{independent}. 
This definition also applies to oriented edges in directed graphs. 
A \emph{clique} in an graph $H$ is a subset $C\subseteq V(H)$ such that all vertices of $C$ are pairwise adjacent. If $H$ is directed, we say that $C$ is a clique (resp. \emph{full clique}) if one of (resp.both) the two possible edges exist between $u$ and $v$, for every pair $u,v\in C$, $u\neq v$. A \emph{biclique} in a bipartite graph $H$ is a disjoint pair of sets $A,B \subseteq V(H)$ such that there is an edge from every $a \in A$ to every $b \in B$. A graph is \emph{$2K_2$-free} if it does not contain a pair of edges $uv$ and $xy$ such that $G[\{u,v,x,y\}]$ contains exactly the two edges (i.e., is isomorphic to a $2K_2$).

A \emph{walk} in a (directed) graph $G$ is an alternating sequence $W$ of vertices and edges that starts and ends with a vertex, and such that for every edge $(u,v)$ in the walk, vertex $u$ (resp. vertex $v$) is the element right before (resp. right after) edge $(u,v)$ in $W$.
If the first vertex in a walk is $u$ and the last one is $v$, then we say this is a \emph{walk from $u$ to $v$}.
A \emph{path} is a (directed) graph containing exactly a walk that contains all of its vertices and edges without repetition. 
A directed graph $D$ is \emph{strongly connected} if, for every pair of vertices $u,v \in V(D)$, there is a walk from $u$ to $v$ and a walk from $v$ to $u$ in $D$.
We say that $D$ is \emph{weakly connected} if the underlying graph of $D$ is connected.
A \emph{strong component} of $D$ is a maximal induced subgraph of $D$ that is strongly connected, and a \emph{weak component} of $D$ is a maximal induced subgraph of $D$ that is weakly connected.
For simplicity, we call a strong connected component of a directed graph simply ``connected component''. 

An \emph{orientation} of an undirected graph $G$ is a digraph $D$ obtained from $G$ by choosing an orientation for each edge $e \in E(G)$.
The undirected graph $G$ formed by ignoring the orientation of the edges of a digraph $D$ is the \emph{underlying graph} of $D$.

The \emph{line graph} of an undirected graph $G$ is the graph $L$ with vertex set $E(G)$ where two vertices $e,f \in V(L)$ associated with edges of $G$ are linked by an edge if and only if $e$ and $f$ share an endpoint in $G$.

For a undirected (static) graph $G$ and $u \in V(G)$, let $\delta(u)$ denote the set $\{e \in E(G) \mid u \in e\}$. 
For a set $S \subseteq E(G)$, the \emph{edge-induced subgraph} $G[S]$ is a graph whose edge set is $S$ and vertex set consists of all endpoints of the edges in $S$.

Given a temporal (directed) graph ${\cal G} = (G,\lambda)$ and a subset $S\subseteq V(G)$, we say that $S$ is \emph{temporal connected} (in ${\cal G}$) if there is a temporal $u,v$-path in ${\cal G}$ for every ordered pair $(u,v)\in S\times S$.

\subsection{Parameterized complexity}

A \emph{parameterized problem} is a language $L \subseteq \Sigma^* \times \mathbb{N}$.  For an instance $I=(x,k) \in \Sigma^* \times \mathbb{N}$, $k$ is called the \emph{parameter}.

A parameterized problem $L$ is \emph{fixed-parameter tractable} ({\sf FPT}) if there exists an algorithm $\mathcal{A}$, a computable function $f$, and a constant $c$ such that given an instance $I=(x,k)$, $\mathcal{A}$   (called an {\sf FPT} \emph{algorithm}) correctly decides whether $I \in L$ in time bounded by $f(k) \cdot |I|^c$. For instance, the \textsc{Vertex Cover} problem parameterized by the size of the solution is {\sf FPT}.

A parameterized problem $L$ is in {\sf XP} if there exists an algorithm $\mathcal{A}$ and two computable functions $f$ and $g$ such that given an instance $I=(x,k)$, $\mathcal{A}$  (called an {\sf XP} \emph{algorithm}) correctly decides whether $I \in L$ in time bounded by $f(k) \cdot |I|^{g(k)}$. For instance,  the \textsc{Clique} problem parameterized by the size of the solution is in  {\sf XP}.

Within parameterized problems, the class {\sf W}[1] may be seen as the parameterized equivalent to the class {\sf NP} of classical decision problems. Without entering into details (see~\cite{DF13,CyganFKLMPPS15} for the formal definitions), a parameterized problem being {\sf W}[1]-\emph{hard} can be seen as a strong evidence that this problem is {\sl not} {\sf FPT}.
\textsc{Clique} parameterized by the size of the solution is the canonical example of a $\W$[1]-hard problem.

\emph{Parameterized reductions} are used to transfer fixed-parameter tractability or hardness between parameterized problems.
Namely, a parameterized reduction is an algorithm  that, given an instance $(x, k)$ of a parameterized problem $L$, runs in time $f(k)\cdot |x|^{O(1)}$ and outputs an instance $(x', k')$ of a parameterized problem $L'$ such that $k' \leq g(k)$ for some computable function $g$ and $(x, k)$ is positive if and only if $(x',k')$ is positive.
For example, if $L$ is \textsf{W}[1]-hard and there is a parameterized reduction from $L$ to $L'$, then $L'$ is also \textsf{W}[1]-hard and thus unlikely to admit an \textsf{FPT} algorithm.

\subsection{Temporal paths and components}
Given a temporal (directed) graph $\mathcal{G} = ( G,\lambda)$ and vertices $v_0,v_{q}\in V(G)$, a \emph{temporal $v_0,v_{q}$-walk} in $\mathcal{G}$ is defined as a sequence of vertices and temporal edges,  $(v_0, \alpha_1,v_1$, $\cdots,\alpha_q,v_q)$ such that, for each $i\in [q]$, $\alpha_i$ has endpoints $v_{i-1}v_{i}$ and is active in a timestep $t_i$ which is at most equal to the timestep where $\alpha_{i+1}$ is active. Sometimes we abuse notation and write $P=(v_0,t_1,v_1,\ldots,t_q,v_q)$ instead, where $t_i$ is equal to the timestep where $\alpha_i$ is active, for every $i\in [q]$. We then say that $P$ \emph{starts in time $t_1$} and \emph{finishes in time $t_q$}. Given $i,j\in \{0,\ldots,q\}$, we denote by $v_iPv_j$ the $v_i,v_j$-walk $(v_i,t_{i+1},v_{i+1},\ldots,t_j,v_j)$. 
Additionally, if no vertices of $G$ are repeated in $P$, then we say that $P$ is a \emph{temporal $v_0,v_q$-path}. It is important to mention that distinctions between paths and walks are important for some problems, but since it is not the case in this work, the reader should not worry about interchangeable uses of ``walks'' and ``paths'' along the text. 
Given two vertices $u,v \in V(G)$, we say that \emph{$u$ reaches $v$} in $\mathcal{G}$ if there exists a temporal $u,v$-walk in $\mathcal{G}$. 

Given a temporal (directed) graph ${\cal G} = (G,\lambda)$ and a subset $S\subseteq V(G)$, we say that $S$ is \emph{temporal connected} (in ${\cal G}$) if there is a temporal $u,v$-path in ${\cal G}$ for every ordered pair $(u,v)\in S\times S$. 
In~\cite{BF.03}, the authors define a \emph{temporal connected component} ({\TCC} for short) as a maximal subset $S\subseteq V(G)$ such that $S$ is temporal connected.
Similarly, a \emph{closed temporal connected component} ({\cTCC} for short) was defined as a maximal subset $S\subseteq V(G)$ for which, for every ordered pair $(u,v)\in S\times S$, there is a temporal $u,v$-path in ${\cal G}$ using only vertices of $S$. In other words, a {\cTCC} is a maximal subset $S\subseteq V(G)$ such that ${\cal G}[S]$ (the temporal subgraph induced by $S$) is temporal connected.
We say that $S$ is a \emph{temporal unilaterally connected set} if for every pair $u,v \in S$ there is a temporal $u,v$-path or a temporal $v,u$-path in $\mathcal{G}$.
If all such paths use only vertices in $S$ then we say that $S$ is a \emph{closed temporal unilaterally connected set}.
If $S$ is maximal such set, then we say that $S$ is a \emph{temporal unilaterally connected component} (\TUCC) in the first case, and a \emph{closed temporal unilaterally connected component} (\cTUCC) in the second case.

The \emph{reachability digraph} \reach$(\mathcal{G})$ associated to $\mathcal{G}$ is a directed graph with the same vertex set as ${\cal G}$, and such that $uv$ is an edge in $\reach(\mathcal{G})$ if and only $u$ reaches $v$ in ${\cal G}$, $u\neq v$. This is related to the \emph{affine} graph in~\cite{nicosia.etal.12}. 
This is a slight generalization of the \emph{affine} graph introduced in~\cite{nicosia.etal.12}. There, since they are interested only in the {\TCC} variants, they consider pairs that are mutually reachable from each other, ignoring the edges $uv$ of $\reach({\cal G})$ that are not symmetric (i.e., for which $vu$ is not present). 

\section{Auxiliary results}

The following result is an immediate consequence of the definition of rechability graph \reach$(\mathcal{G})$.

\begin{lemma}\label{lemma:cliquecomponent}
Given a temporal (directed) graph $\mathcal{G} = ( G,\lambda)$, then the following hold: 
\begin{enumerate}
    \item $C$ is a {\TCC} in ${\cal G}$ if and only if $C$ is a maximal full clique in \reach$(\mathcal{G})$; 
    \item $C$ is a {\TUCC} in ${\cal G}$ if and only if $C$ is a maximal clique in \reach$(\mathcal{G})$;
    \item  If $C$ is a {\cTCC} in ${\cal G}$, then $C$ is a full clique in \reach$(\mathcal{G})$; and
    \item If $C$ is a {\cTUCC} in ${\cal G}$, then $C$ is a clique in \reach$(\mathcal{G})$.
\end{enumerate}
\end{lemma}

For each $i\in [\tau]$ and $u\in V({\cal G})$, denote by $C_i(u)$ the set of vertices in the same connected component of $G_i$ as $u$, and by $R_i(u)$ the set of vertices in $G_i$ reachable from $u$ (i.e., $v\in R_i(u)$ if and only if there is a $u,v$-path in $G_i$). Observe that, if $G$ is undirected, then $C_i(u) = R_i(u)$. 
Note also that ${\cal R}_\tau(u)$ is exactly equal to $N^-_{\reach({\cal G})}(u)$.
For the sake of completeness, we now show that we can recursively define the set of vertices reachable from $u$ by a temporal path finishing at time at most $i$. We apply the following lemma in the context of non-strict reachability, but also holds for strict. 

\begin{lemma}\label{lem:reachableset}
Let $\mathcal{G}$ be a (directed) temporal graph, and let ${\cal R}_i(u)$ be recursively defined as:
\[
{\cal R}_i(u) = \left\{
\begin{array}{ll}
      R_1(u) & \mbox{, if $i=1$}\\ 
      \bigcup_{v\in {\cal R}_{i-1}(u)} R_i(v) & \mbox{, otherwise}
\end{array}\right.\]
Then ${\cal R}_i(u)$ is equal to the set of vertices reachable from $u$ by a temporal path finishing at time at most $i$. 
\end{lemma}
\begin{proof}
We want to prove that $v\in {\cal R}_i(u)$ if and only if there exists a temporal $u,v$-walk finishing in time at most $i$. 
First, let $v\in {\cal R}_i(u)$. If $i = 1$, then $v \in R_1(u)$ and $u$ reaches $v$ in $G_1$ by definition. So suppose $i > 1$. 
Again by definition, we have $v\in \bigcup_{v'\in {\cal R}_{i-1}(u)} R_i (v')$. 
Consider then $w \in {\cal R}_{i-1}(u)$ such that $v \in R_{i}(w)$. By induction hypothesis, there exists a temporal $u,w$-path $P$ finishing in time at most $i-1$. And because $w$ reaches $v$ in $G_i$, such path can be extended to a temporal $u,v$-walk finishing in time at most $i$. 

Now, let $v$ be such that there exists a temporal $u,v$-path $P$ finishing in time at most $i$. If $P$ finishes in time at most $i-1$, we are done by induction hypothesis. Otherwise, let $w\in V(P)$ be closest to $v$ in $P$ such that the temporal edges incident to $w$ in $P$ occur in time $j<i$ and $i$. Observe that $wPv$ is contained in $G_i$, and hence $v\in R_i(w)$. Additionally, $uPw$ finishes in time at most $j\le i-1$, and by induction hypothesis $w\in {\cal R}_{i-1}(u)$. By definition we then get $v\in {\cal R}_i(u)$, as we wanted to show. 
\end{proof}

The following easy proposition tells us that deciding the existence of large components (i.e. maximal connected sets) is equivalent to deciding the existence of large connected sets.

\begin{proposition}\label{prop:components_vs_sets}
    Let ${\cal G}$ be a temporal (directed) graph. Then ${\cal G}$ has a (closed) temporal (unilaterally) connected component of size at least $k$ if and only if (closed) temporal (unilaterally) connected set of size at least $k$.
\end{proposition}


\section{Parameterized Complexity Results: Proofs}

\subsection{Proof of Theorem~\ref{thm:param_tau}}\label{app:param_tau}

The \textsc{Maximum Edge Biclique} Problem (\textsc{MEBP} for short) consists in, given a bipartite graph $G$ and an integer $k$, deciding whether $G$ has a biclique with at least $k$ edges. It was proved to be {\NP}-complete in~\cite{P.03}.
Problem \textsc{$2K_2$-free Edge Subgraph} consists in, given a bipartite graph $G$ and an integer $k$, deciding whether $G$ has a $2K_2$-free subgraph with at least $k$ edges. 
This was proved to be $\NP$-complete in~\cite{Y.81}.

For an undirected (static) graph $G$ and $u \in V(G)$, let $\delta(u)$ denote the set $\{e \in E(G) \mid u \in e\}$. And for a set $S \subseteq E(G)$, the \emph{edge-induced subgraph} $G[S]$ is a graph whose edge set is $S$ and vertex set consists of all endpoints of the edges in $S$.

As stated in the main text, we prove {\NP}-completeness of {\gTCC} and {\gTUCC} using the same construction and reducing from the above problems. We start by proving hardness of {\gTCC}.

\begin{proof}
Consider an instance $(H,k)$ of \textsc{MEBP}, consisting of a bipartite graph $H=(X \cup Y,E)$ and an integer $k$. Let $X=\{x_1,\ldots,x_p\}$ and $Y=\{y_1,\ldots,y_q\}$. We construct a temporal graph $\mathcal{G}=(G,\lambda)$ with lifetime~2 such that $V(G)=E(H)$ and snapshot $G_1$ is the graph whose connected components are precisely $\delta(x_i)$ for each $i \in \{1,\ldots,p\}$, while snapshot $G_2$ is the graph whose connected components are precisely $\delta(y_j)$ for each $j \in \{1,\ldots,q\}$. We consider that the components are cliques; clearly $G$ is the line graph of $H$. 
We claim that there exists a biclique $(A,B)$ in $H$ with at least $k$ edges if and only if there exists a {\cTCC} in $\mathcal{G}$ of size at least $k$. Then we prove that every {\TCC} is also a {\cTCC}, finishing this part of the proof.

Suppose first that there exists a biclique $(A,B)$ in $H$ with at least $k$ edges, and let $C=E(H[A,B])$. 
We want to show that $C$ is a closed temporal connected set of $\mathcal{G}$. Let $e = xy$ and $e' = x'y'$ be two elements of $C$, with $\{x,x'\}\subseteq X$ and $\{y,y'\}\subseteq Y$. 
If $x=x'$ then $e,e'\in \delta(x)$ and hence are contained in the same component of $G_1$ (i.e., they reach each other by a direct edge); the analogous holds in case $y=y'$, so suppose $x\neq x'$ and $y\neq y'$. 
Since $(A,B)$ is a biclique in $H$, we have that $\{xy,xy',x'y,x'y'\} \subseteq C$. Denote $xy'$ by $f$ and $x'y$ by $f'$. 
Now, in $\mathcal{G}$ we can reach $f$ from $e$ at timestep $1$ and $e'$ from $f$ at timestep $2$.
Similarly, we can also reach $e$ from $e'$ in $\mathcal{G}$.
Because $f,f'$ are also in $C$, and since this holds for any two such edges, we get that $C$ is a closed temporal connected set, and by Proposition~\ref{prop:components_vs_sets}, we get that $\mathcal{G}$ has a {\cTCC} of size at least $k$. 

For the converse, suppose that $\mathcal{G}$ has a {\TCC} $C$ with $|C| \ge k$. We want to show that $C$ forms a biclique in $H$ with at least $k$ edges. Let $A\subseteq X$ contain all vertices incident to some $e\in C$, and define $B$ similarly with relation to $Y$. 
First we show that $(A,B)$ is a biclique in $H$. Observe that, combined with the previous paragraph, we get that any {\TCC} is also a {\cTCC}; hence the proof will follow also for both problems. Let $x \in A$ and $y\in B$. We need to show that $xy$ is an edge of $H$. 
Note that since $x \in A$, it must be an endpoint of some edge $e_x \in C$; analogously, since $y \in B$, it must be an endpoint of some edge $e_y \in C$. Let $y'$ be the other endpoint of $e_x$ and let $x'$ be the other endpoint of $e_y$. 
Since $H$ is bipartite, we have $x' \in A$ and $y' \in B$. If $x'=x$ or $y'=y$ we are done, so suppose otherwise. Since $e_x$ and $e_y$ are in $C$, there exists a temporal $e_x,e_y$-path. Note that, by the construction of $\mathcal{G}$, this means that there exists an edge $e' \in \delta(x) \cap \delta(y)$. The only possibility is $e' = xy$, as we wanted. To finish the proof just observe that $\lvert E(H[A,B])\rvert \ge \lvert C\rvert\ge k$.

For {\gTUCC}, we make a reduction from \textsc{$2K_2$-free Edge Subgraph}. 
This was shown to be equivalent to the \textsc{Minimum Fill-in} problem in co-bipartite graphs in~\cite{Y.81}, where the authors also showed that this problem is $\NP$-complete. 
Given a bipartite $H$, the proof follows similarly to the first case, using exactly the same construction for the temporal graph $\mathcal{G}$.

Assume that $H$ contains a $2K_2$-free subgraph $H'$ with at least $k$ edges.
Given edges $e = xy$ and $e' = x'y'$ of $H'$, since $H'$ is $2K_2$-free either $xy' \in E(H')$, or $x'y \in E(H')$, or both. 
Thus in $\mathcal{G}$ there is a temporal path from $e$ to $e'$, or a temporal path from $e$ to $e'$, or both.
Since in this case we only care about unilateral components, the first part of the proof follows.
The second part also follows with similar arguments, just noticing that an unilateral component only generates a $2K_2$-free subgraph in $H$.

In both cases, for higher values of $\tau$ is suffices to add snapshots with empty edge sets.
\end{proof}


\subsection{Proof of Theorem~\ref{thm:TCCbyk}}\label{app:TCC_Undirected_Whard}

\begin{proof}
We make a parameterized reduction from $k$-\textsc{Clique}.
Let $G$ be graph and $k \geq 3$ be an integer.
We construct the temporal graph $\mathcal{G} = (G', \lambda)$ as follows.
Fist, add to $G'$ every vertex in $V(G)$ and make $V = V(G)$.
Second, add to $G'$ a copy $u'$ of every vertex $u \in V$ and define $V' = \{u' \mid u \in V\}$.
Third, for every pair $u,u'$ with $u\in V$ and $u' \in V'$ add the edge $uu'$ to $G'$ and make all such edges active at timestep $0$.
Fourth, consider an arbitrary ordering $e_1, \ldots, e_m$ the edges of $G$ and
for each edge $e_i = uv$ create for new vertices $\{h_{uv},h_{vu},h'_{uv},h'_{vu}\mid uv\in E(G)\}$, and add edges:
 \begin{itemize}
  \item $uh_{uv}$ and $vh_{vu}$, active at time $i$;
  \item $u'h'_{uv}$ and $v'h'_{vu}$, active at time $2m+i$;
  \item $h_{vu}u$ and $h_{uv}v$, active at time $m+i$; and
  \item $h'_{vu}u'$ and $h'_{uv}v'$, active at time $3m+i$.
\end{itemize}
Denote the set $\{h_{uv},h_{vu}\mid uv\in E(G)\}$ by $H$, and the set $\{h'_{uv},h'_{vu}\mid uv\in E(G)\}$  by $H'$.
We now prove that $G$ has a clique of size at least $k$ if and only if ${\cal G}$ has a temporal connected set of size at least $2k$. The theorem follows by Proposition~\ref{prop:components_vs_sets}.

\begin{figure}[t]
\begin{center}
\scalebox{.6}{
\begin{tikzpicture}[scale=1,vertex/.style={circle, minimum size=0.2cm, draw, inner sep=1pt}, blackvertex/.style={draw,circle,minimum size=5pt,inner sep=0pt, fill=black}]

\begin{scope}[xshift = -3cm]
\draw[rounded corners] (0,0) rectangle  (1.2,5) node [above, xshift = -.5cm] {\large$V$};

\draw[rounded corners] (2,0) rectangle  (3.2,5) node [above, xshift = -.5cm] {\large$V'$};

\foreach \x/\i in {0.5/n, 2.5/3, 3.5/2, 4.5/1}{
\node[blackvertex, label = left:{\large$u_{\i}$}] (u\i) at (0.75, \x) {};
\node[blackvertex, label = right:{\large$u'_{\i}$}] (v\i) at (2.45, \x) {};
\draw[thick] (u\i) -- (v\i) node[midway, fill=white] {$0$};
}

\path (u3) -- (un) node [font=\Huge, midway, sloped] {$\dots$};
\path (v3) -- (vn) node [font=\Huge, midway, sloped] {$\dots$};
\end{scope}

\begin{scope}[xshift=6cm, yshift = 5.5cm]
\node[blackvertex, scale = 1.1, label={\large$u$}] (u) at (0,0) {};
\node[blackvertex, scale = 1.1, label={\large$v$}] (v) at (2,0) {};
\draw[thick] (u) -- (v) node[midway, above]{\large$e_i$};
\end{scope}

\begin{scope}[xshift = 6cm, yshift = 3.5cm]
\node[blackvertex, scale = 1.1, label=left:{\large$u$}] (u) at (-2,0) {};
\node[blackvertex, scale = 1.1, label=right:{\large$v$}] (v) at (4,0) {};
\node[blackvertex, scale = 1.1, label=below:{\large$h_{uv}$}] (huv) at (1,1) {};
\node[blackvertex, scale = 1.1, label={\large$h_{vu}$}] (hvu) at (1,-1) {};

\draw[thick] (u) -- (huv) node[midway, fill=white] {$i$};
\draw[thick] (hvu) -- (v) node[midway, fill=white] {$i$};
\draw[thick] (huv) -- (v) node[midway, fill=white] {$m+i$};
\draw[thick] (u) -- (hvu) node[midway, fill=white] {$m+i$};

\begin{scope}[yshift = -2.5cm]
\node[blackvertex, scale = 1.1, label=left:{\large$u'$}] (u1) at (-2,0) {};
\node[blackvertex, scale = 1.1, label=right:{\large$v'$}] (v1) at (4,0) {};
\node[blackvertex, scale = 1.1, label=below:{\large$h'_{uv}$}] (huv1) at (1,1) {};
\node[blackvertex, scale = 1.1, label={\large$h'_{vu}$}] (hvu1) at (1,-1) {};

\draw[thick] (u1) -- (huv1) node[midway, fill=white] {$2m+i$};
\draw[thick] (hvu1) -- (v1) node[midway, fill=white] {$2m+i$};
\draw[thick] (huv1) -- (v1) node[midway, fill=white] {$3m+i$};
\draw[thick] (u1) -- (hvu1) node[midway, fill=white] {$3m+i$};
\end{scope}
\draw[thick] (u) -- (u1) node [midway, fill=white] {$0$};
\draw[thick] (v) -- (v1) node [midway, fill=white] {$0$};

\end{scope}

\end{tikzpicture}
}
\end{center}
\caption{This appears as Figure~\ref{fig:theorem-2} in the main text. Construction used in the proof of Theorem~\ref{thm:TCCbyk}. On the left, the two copies of $V(G)$ and the edges between them, active in timestep $0$. On the right, the edge $e_i \in E(G)$ and the associated gadget in $\mathcal{G}$.}
\label{fig_app:theorem-2}
\end{figure}
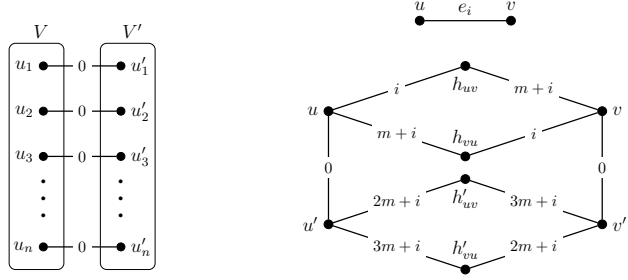

First, let $C\subseteq V$ be a clique of size at least $k$ in $G$ and $C' = \{u'\mid u\in C\}$. We show that $C \cup C'$ is a \TCC of $\mathcal{G}$.
For this, let $u,v\in C$. Since $\lambda(uh_{uv}) = i < m+i = \lambda(h_{uv}v)$, we get that $u$ reaches $v$ in $\mathcal{G}$ through $h_{uv}$.
Because $C$ is a clique in $G$ we conclude that $C$ is a temporal connected set of $\mathcal{G}$, and similarly the same holds for $C'$.
Thus it remains to show that pairs of vertices of the form $u,u'$ with $u \in C$ and $u'\in C'$ are also connected in $\mathcal{G}$.
This is true due to the choice of timestep $0$ for the edges forming the matching between $V$ and $V'$ of $\mathcal{G}$.

Now, let $S\subseteq V(G')$ be a \TCC of ${\cal G}$ of size at least $2k$. 
We want to show that either $C = \{u\in V(G)\mid u\in S\cap V\}$ or $C'= \{u\in V(G)\mid u'\in S\cap V'\}$ is a clique of $G$ of size at least $k$. 
For this, we first prove a series of useful facts.

\begin{clm}\label{claim:paths_size}
Let $P$ be a temporal path in $G'[V\cup H]$.
Then $P$ has at most one internal vertex of $H$, and hence $\lvert V(P)\rvert\le 5$. The same holds if $P$ is contained in $G'[V'\cup H']$.
\end{clm}
\begin{proof}
It suffices to observe that every $e\in H$ is incident to exactly two edges of $G'$, one active at time at most $m$ and the other one active at time at least $m+1$. The second part follows because $G'[V\cup P]$ is a bipartite graph.
\end{proof}

\begin{clm}\label{claim:Visclique}
$C$ and $C'$ are cliques in $G$.
\end{clm}
\begin{proof}
Let $u,v\in C$.
Since $C$ is a temporal connected set, there is a temporal path from $u$ to $v$.
Such path must contain only edges of $G'[V\cup H]$ since the edges between $V$ and $V'$ are only active in timestep $0$, and all other edges are active in a later time (i.e., there is no way to leave $u$ to $u'$ at time $0$, then go back to $v$). 
By Claim~\ref{claim:paths_size} and the fact that $G'[V\cup H]$ is bipartite, it follows that $u$ and $v$ must be adjacent. 
The argument for $u,v\in C'$ is analogous by taking their copies, $u',v'$ in $S$. 
\end{proof}

Note that if $S\subseteq V\cup V'$, then  Claim~\ref{claim:Visclique} and the fact that $\lvert S\rvert \ge 2k$ directly imply that either $C$ or $C'$ is a clique of size at least $k$ in $G$. 
Assume now that $S\cap (H\cup H')\neq \emptyset$.
In this case, it is not ensured that $C$ or $C'$ contains a clique of size at least $k$, but the following claims allow us to obtain another clique.

\begin{clm}\label{claim:H_to_VHprime}
For every $h\in H$ and every $x'\in V'\cup H'$, $h$ does not reach $x'$. Similarly, for every $h'\in H'$ and every $x\in V\cup H$, $h'$ does not reach $x$.
\end{clm}
\begin{proof}
The only edges between $V\cup H$ and $V'\cup H'$ are those incident to $V$ and $V'$ at timestep $0$.
Since every edge incident to $h\in H\cup H'$ is active only at a later timestep, the claim follows.
\end{proof}

\begin{clm}\label{claim:intersectionwithH}
If $a,b\in S\cap H$, then $a$ and $b$ are related to the same edge, or to edges adjacent to each other. The same holds for $a,b\in S\cap H'$.
\end{clm}
\begin{proof}
Suppose, without loss of generality, that $a$ reaches $b$.
Suppose also by contradiction that $a,b$ are related to distinct edges, say $e_i$ and $e_j$, respectively.
Write $e_i$ as $uv$ and $e_j$ as $xy$ and assume that that $\{u,v\} \cap \{x,y\} = \emptyset$.
Because every temporal path between $a$ and $b$ must alternate between $V$ and $H$, as $G'[V\cup H]$ is bipartite, and since by Claim~\ref{claim:paths_size} every temporal path contains at most one internal vertex of $H$, we get that the temporal $a,b$-path must use vertices $(a,v,h_{vx},x,b)$. 
This gives us that $ev$ and $vh_{vx}$ must be active in timestep at most $m$, while $h_{vx}x$ and $xb$ must be active in timestep at least $m+1$. 
Hence, by letting $vx$ be equal to $e_\ell$, we must have that $i < \ell < j$.
We apply an analogous argument to a temporal $b,a$-path to obtain that $j$ must be smaller than $i$, a contradiction.
A similar argument can clearly  be applied to $e,f\in S\cap H'$, and the claim follows.
\end{proof}

Now suppose that $S\cap H\neq \emptyset$. By Claim~\ref{claim:H_to_VHprime} we get that $S\subseteq V\cup H$. Since $V$ and $H$ are disjoint and $\lvert S\rvert \ge 2k$, we get that either $\lvert S\cap V\rvert\ge k$ or $\lvert S\cap H\rvert\ge k$. If the former occurs, then $C$ contains a clique of size at least $k$ by Claim~\ref{claim:Visclique}. Otherwise, denote by $E_S$ the set of edges of $G$ related to vertices in $S\cap H$ (i.e. $E_S = \{uv\in E(G)\mid \{h_{uv},h_{vu}\}\cap S\neq \emptyset\}$). The following is the last ingredient of the proof.

\begin{clm}\label{claim:triangles}
Let $a,b \in S \cap H$ be associated with distinct edges $g,g'$ of $G$ sharing an endpoint $v$.
If $u$ and $w$ are the other endpoints of $g$ and $g'$, respectively, then $u$ and $w$ are also adjacent in $G$.
Additionally, either $\lvert S\cap \{h_{xy},h_{yx}\}\rvert \le 1$ for every $xy\in E(G)$, or $\lvert S\cap H\rvert \le 2$.
\end{clm}
\begin{proof}
Suppose, without loss of generality, that $a$ reaches $b$.
By contradiction suppose that $u$ and $w$ are not adjacent in $G$.
This gives us that every $uw$-path in $G'$ contains two internal vertices of $H$, and therefore is not a temporal path by Claim~\ref{claim:paths_size}. Because every subpath of a temporal path is also a temporal path, this means that there is no temporal $a,b$-path passing by $u$ and $w$. 
By construction, and since $G$ is a simple graph (i.e., there is only one edge with endpoints $u$ and $v$, and only one with endpoints $v$ and $w$), we get from Claim~\ref{claim:paths_size} that, if $P$ is a temporal $a,b$-path not containing both $u$ and $w$, then $P$ is one of the following paths: $P_1 = (a,v,b)$; $P_2 = (a,u,a',v,b)$ where $\{a,a'\} = \{h_{uv},h_{vu}\}$; or $P_3 = (a,v,f',w,f)$  where $\{f,f'\} = \{h_{vw},h_{wv}\}$. 
Note that the same argument can be applied to a temporal $b,a$-path, except that we arrive to the reverses of the above paths. 
Note that, since all paths are strictly increasing, we get that at least one between $P_2$ or $P_3$ (or their reverse) is a temporal path. But observe that neither the subpath  $(a,u,a',v)$ nor its reverse can ever be temporal paths by construction, which means that neither $P_2$ nor its reverse can be temporal paths. 
A similar argument can be applied to $P_3$, thus leading to a contradiction.

For the second part, suppose by contradiction that $\{h_{xy},h_{yx}\}\subseteq S$ and $\lvert S\cap H\rvert > 2$.
Let $a \in (S\cap H)\setminus \{h_{xy},h_{yx}\}$.
By Claim~\ref{claim:intersectionwithH} we can suppose, without loss of generality, that $a \in \{h_{xw},h_{wx}\}$ for some $w \neq y$.
Observe also that the previous paragraph tells us that one of the temporal paths between $\{h_{xy},h_{yx}\}$ and $a$ must contain $(y,f,w)$ or its reverse, where $f\in \{h_{yw},h_{wy}\}$.
Since such a path contains $4$ edges, by letting $xy$ be equal to $e_i$, $yw$ be equal to $e_j$ and $wx$ be equal to $e_\ell$, we get $i<j<\ell$. 
Thus in this case we have that $we$ is active in time at least $m+1$, which in turn gives us that $a = h_{xw}$.
We can now verify that $a$ does not reach $h_{xy}$.
Indeed, every $a,h_{xy}$-path starting with edge $aw$ must contain some internal vertex $h$ of $H$, in which case it cannot be a temporal path as it starts with an edge active at time at least $m+1$ (namely $aw$) and contains an edge active in time at most $m$ (namely one of the edges incident to $h$). A similar argument can be applied if the path starts with edge $ax$, since it must be distinct from $(e,x,h_{xy})$ (recall that $\lambda(ax) = \ell > i = \lambda(xh_{xy})$).
\end{proof}

Now, recall that we are in the case $\lvert S\cap H\rvert \ge k+1$. By our assumption that $k\ge 3$, note that Claim~\ref{claim:triangles} gives us that $\lvert S\cap \{h_{xy},h_{yx}\}\rvert \le 1$ for every $xy\in E(G)$, which in turn implies that $\lvert E_S\rvert = \lvert S\cap H\rvert$. Additionally, observe that, since $\lvert S\cap H\rvert \ge 4$, Claim~\ref{claim:triangles} also gives us that there must exist $w\in V$ such that $e$ is incident to $w$ for every $e\in E_S$. 
Indeed, the only way that $3$ distinct edges can be mutually adjacent without being all incident to a same vertex is if they form a triangle.
Supposing that $3$ edges in $E_S$ form a triangle $T = (a,b,c)$, since $\lvert E_S\rvert \ge 4$, there exists an edge $e\in E_S\setminus E(T)$.
But now, since $G$ is a simple graph, $e$ is incident to at most one between $a$, $b$ and $c$, say $a$.
We get a contradiction to Claim~\ref{claim:triangles} as in this case $e$ is not incident to edge $bc\in E_S$. 
Finally, by letting $C'' = \{v_1,\ldots,v_k\}$ be any choice of $k$ distinct vertices such that $\{wv_1,\ldots,wv_k\}\subseteq E_S$, Claim~\ref{claim:triangles} gives us that $v_i$ and $v_j$ are adjacent in $G$, for every $i,j\in [k]$; i.e., $C''$ is a clique of size at least $k$ in $G$.
This finishes the proof as the case $S\cap H'\neq \emptyset$ is clearly analogous. 
\end{proof}


\subsection{Proof of Theorem~\ref{thm:W1h_dir_tau2}}
\label{app:W1h_dir_tau2}

\begin{proof}
See Figure~\ref{fig_app:thm_3} to follow the construction. Let $G$ be a graph and consider the directed graph $D_G$ constructed as follows. 
First, add to $D_G$ every vertex of $G$. 
Then, for each $uv \in E(G)$, add to $D_G$ vertices $h_{uv}$ and $h_{vu}$, directed edges $uh_{uv}$ and $vh_{vu}$, and directed edges $h_{uv}v$ and $h_{vu}u$. 
Denote by $H$ the set $\{h_{uv},h_{vu}\mid uv\in E(G)\}$. 
To construct the directed temporal graph $\mathcal{G}$ we start from $D_G$ and for every $uv \in E(G)$
\begin{itemize}
    \item make edges $uh_{uv}$ and $vh_{vu}$ active in timestep $1$; and
    \item make edges $h_{uv}v$ and $h_{vu}u$ active in timestep $2$.
\end{itemize}

\begin{figure}[t]
    \centering
    \begin{subfigure}[b]{0.4\textwidth} {\begin{tikzpicture}[scale=1]

\node (G) at (-.75,1.75) {$G$};

\vertex[label=left:$u$] (u) at (0,0) {};
\vertex[label=left:$v$] (v) at (0,1) {};
\vertex[label=left:$z$] (z) at (0,2) {};

\draw (u) -- (v) -- (z);

\end{tikzpicture}}
    \caption{}\label{fig_app:graphG}
    \end{subfigure}
    \begin{subfigure}[b]{0.4\textwidth}{\begin{tikzpicture}[scale=1]
 \node (G) at (-.5,-.75) {$\mathcal{G}$};
\vertex[label=left:$u$] (u) at (0,0) {};
\vertex[label=above:$h_{uv}$] (huv) at (.75,.75) {};
\vertex[label=below:$h_{vu}$] (hvu) at (.75,-.75) {};
\vertex[label=right:$v$] (v) at (1.5,0) {};
\vertex[label=above:$h_{vz}$] (hvz) at (2.25,.75) {};
\vertex[label=below:$h_{zv}$] (hzv) at (2.25,-.75) {};
\vertex[label=right:$z$] (z) at (3,0) {};

 \draw[->] (u) to node[scale=1,sloped,above] {$1$} (huv);
 \draw[->] (huv) to node[scale=1,sloped,below] {$2$}  (v);
 \draw[->] (v) to node[scale=1,sloped,below] {$1$} (hvu);
 \draw[->] (hvu) to node[scale=1,sloped,below] {$2$} (u);
 \draw[->] (v) to node[scale=1,sloped,above] {$1$} (hvz);
 \draw[->] (hvz) to node[scale=1,sloped,above] {$2$} (z);
 \draw[->] (z) to node[scale=1,sloped,below] {$1$} (hzv);
 \draw[->] (hzv) to node[scale=1,sloped,below,pos=.3] {$2$} (v);

\end{tikzpicture}}
    \caption{}
    \end{subfigure}
    \caption{Given the graph in (a), Theorem~\ref{thm:W1h_dir_tau2} constructs the directed temporal graph in (b).}
    \label{fig_app:thm_3}
\end{figure}
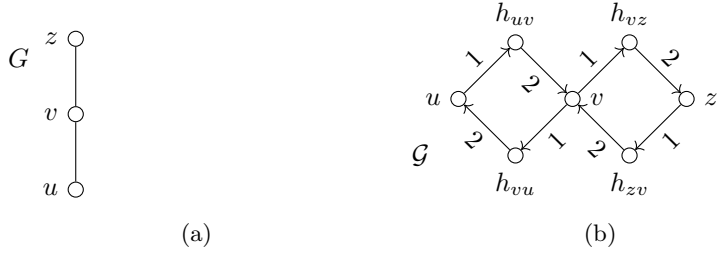

Assume $k\ge 3$. We now prove that $G$ has a clique of size $k$ if and only if $\mathcal{G}$ has a temporal connected set of size at least $k$. The theorem follows by Proposition~\ref{prop:components_vs_sets}. 
Notice that every vertex of $G$ is contained in $V(\mathcal{G})$, and that $\mathcal{G}$ has lifetime~$2$.

If $C$ is a clique in $G$, then for every $u,v\in C$, we get that $u$ reaches $v$ and $v$ reaches $u$ in ${\cal G}$ because of the paths $(u,1,h_{uv},2,v)$ and $(v,1,h_{vu},2,u)$.
It remains to show that if ${\cal G}$ has a temporal connected set of size at least $k$, then $G$ has a clique of size at least $k$.
Let $C'$ be such a temporal connected set.
We prove that $C' \subseteq V(G)$ and $uv\in E(G)$ for every $u,v\in C'$.
First observe that $G_1$ has only edges from $V(G)$ to $H$, and $G_2$, from $H$ to $V(G)$.
This implies that a temporal path must be of length at most $2$, which in turn implies that $C' \cap V(G)$ must be a clique.
Now suppose that there exists $h_{uv}\in C' \cap H$. 
Observe that $h_{uv}$ has exactly one incoming edge, active in timestep $1$, and exactly one outgoing edge, active in timestep $2$. Additionally, observe that every edge outgoing from $v$ is active in timestep $1$.
This means that $v$ is the only vertex of $V({\cal G})$ reachable from $h_{uv}$, contradicting the fact that $k \ge 3$.
Thus we conclude that $C'\subseteq V(G)$ and the result follows.

Now, for the unilateral case, observe that every \TCC is also a \TUCC, hence from the above paragraph we get that if $G$ has a clique of size at least $k$, then $\mathcal{G}$ has a \TUCC of size at least $k$. Now, if $\mathcal{G}$ has a \TUCC of size at least $k$, then observe that the same arguments as before can be applied. Indeed, if $u,v\in C'\cap V(G)$, then it must be that either $u$ reaches $v$ or $v$ reaches $u$, and in any case we have $uv\in E(G)$. Additionally, we know that $C'$ cannot contain any vertex of $H$, as $k\ge 3$ and $v$ is the only vertex reachable by $h_{uv}$ for every $h_{uv}\in H$. 
\end{proof}


\subsection{Proof of Theorem~\ref{thm:cTCC_cTUCC_Whard}}
\label{app:cTCC_cTUCC_Whard}

\begin{proof}
Observe Figure~\ref{fig_app:thm_4} to follow the construction. Let $G$ be a graph and consider the directed graph $D_G$ constructed as follows.
For every $u \in V(G)$, add to $D_G$ vertices $u^{\text{in}}$ and $u^{\text{out}}$, an edge from $u^{\text{in}}$ to $u^{\text{out}}$, and an edge from $u^{\text{out}}$ to $u^{\text{in}}$ (notice that each pair $u^{\text{in}},u^{\text{out}}$ induce a cycle in $\mathcal{G}$). 
Then, for each edge $uv \in E(G)$, add to $D_G$ an edge from $u^{\text{out}}$ to $v^{\text{in}}$ and an edge from $v^{\text{out}}$ to $u^{\text{in}}$.
The directed temporal graph $\mathcal{G}=(D_G,\lambda)$ is such that $\lambda$ is defined as follows. 
\begin{itemize}
    \item For every $u \in V(G)$, make edges between $u^{\text{in}}$ and $u^{\text{out}}$ active in timesteps $1$ and $3$ in both directions; and
    \item For every $uv \in E(G)$, make the edges from $u^{\text{out}}$ to $v^{\text{in}}$ and from $v^{\text{out}}$ to $u^{\text{in}}$ active in timestep $2$.
\end{itemize}

\begin{figure}[t]
    \centering
    \scalebox{1.3}{\begin{tikzpicture}[scale=1]
 \node (G) at (-.5,-.75) {$\mathcal{G}$};
\vertex[label=below:$u^{in}$] (u-) at (0,0) {};
\vertex[label=above:$u^{out}$] (u+) at (0,1) {};
\vertex[label=above:$v^{in}$] (v-) at (1,1) {};
\vertex[label=below:$v^{out}$] (v+) at (1,0) {};
\vertex[label=below:$z^{in}$] (z-) at (2,0) {};
\vertex[label=above:$z^{out}$] (z+) at (2,1) {};

  \draw[->] (u-) to node[scale=1,sloped,anchor=north,below,pos=.5] {$1,3$} (u+);
  \draw[->] (u+) to[out=240,in=120] node[scale=1,sloped,anchor=east,below,pos=.5] {$1,3$}  (u-);
  \draw[->] (u+) to[out=30,in=150] node[scale=1,sloped,anchor=east,above,pos=.5] {$2$} (v-);
   \draw[->] (v+) to node[scale=1,sloped,anchor=west,below,pos=.55] {$\!2$} (u-);
  \draw[->] (v+) to node[scale=1,sloped,anchor=west,below,pos=.5] {$1,3$} (v-); 
  \draw[->] (v-) to[out=240,in=120] node[scale=1,sloped,below,pos=.5] {$~1,3$}  (v+);
   \draw[->] (z+) to[out=150,in=30] node[scale=1,sloped,anchor=east,above,pos=.5] {$2$} (v-);
  \draw[->] (v+) to node[scale=1,sloped,anchor=east,below,pos=.5] {$2$} (z-);
  \draw[->] (z-) to node[scale=1,sloped,anchor=west,above,pos=.5] {$1,3$} (z+);
  \draw[->] (z+) to[out=-60,in=60] node[scale=1,sloped,anchor=east,above,pos=.5] {$1,3$}  (z-);

\end{tikzpicture}}
    \caption{Temporal graph constructed in the Proof of Theorem~\ref{thm:cTCC_cTUCC_Whard}, given the graph given in Figure~\ref{fig_app:graphG}.}
    \label{fig_app:thm_4}
\end{figure}
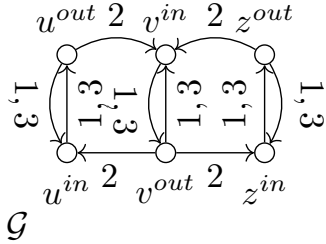

We now prove that $G$ has a clique of size at least $k$ if and only if $\mathcal{G}$ has a closed temporal connected set of size at least $2k$. The theorem follows by Proposition~\ref{prop:components_vs_sets}. Notice that $\mathcal{G}$ has lifetime~$3$. Let $C$ be a clique of size at least $k$ in $G$, and let $C' = \{u^{\text{in}}, u^{\text{out}} \in V(\mathcal{G}) \mid u \in C\}$.
We prove that, for every $u,v\in C$ with $u\neq v$, the set $\{u^{\text{in}}, u^{\text{out}},v^{\text{in}}, v^{\text{out}}\}$ is a closed temporal connected set; note that this implies that $C'$ itself is a closed temporal connected set, as desired. 
By construction, for every $w \in V(G)$ there are temporal paths from $w^{\text{in}}$ to $w^{\text{out}}$ and the other way around, in other words $u^{\text{in}}$ reaches $u^{\text{out}}$, and vice-versa, and $v^{\text{in}}$ reaches $v^{\text{out}}$ and vice-versa. 
Moreover, $u^{\text{in}}$ reaches $v^{\text{in}}$ in $\mathcal{G}$ through the path $(u^{\text{in}},1,u^{\text{out}},2,v^{\text{in}})$. Observe that this also implies that $u^{\text{out}}$ reaches $v^{\text{in}}$, and by symmetry, that both $v^{\text{in}}$ and $v^{\text{out}}$ reach $u^{\text{in}}$. 
Finally, note that the path $(u^{\text{in}},1,u^{\text{out}},2,v^{\text{in}},3,v^{\text{out}})$ implies that both $u^{\text{in}}$ and $u^{\text{out}}$ reach $v^{\text{out}}$, and by symmetry we also get that $v^{\text{in}}$ and $v^{\text{out}}$ reach $u^{\text{out}}$. This finishes this part of the proof.

Assume now that that $C'$ is a closed temporal connected set of ${\cal G}'$ of size at least $2k$.
Let $C = \{u \in V(G) \mid \{u^{\text{in}}, u^{\text{out}}\}\cap C'\neq \emptyset\}$. 
Clearly $|C| \geq k$ since $|C'| \geq 2k$.
To show that $C$ is a clique in $G$, observe that $\mathcal{G}$ consists of a matching at times $1$ and $3$, containing only edges of the form $u^{\text{in}}u^{\text{out}}$ and of the form $u^{\text{out}}u^{\text{in}}$, together with edges in timestep~2 that go only from $O = \{u^{\text{out}}\mid u\in V(G)\}$ to $I = \{u^{\text{in}}\mid u\in V(G)\}$. 
This implies that any temporal path in $\mathcal{G}$ contains at most one edge from $O$ to $I$, which are only defined if the corresponding vertices are adjacent in $G$. We then get that, if $u,v\in C$ with $u\neq v$, then it must be the case that $uv\in E(G)$. 

The proof for {\cTUCC} is similar, except that, for every $uv\in E(G)$, we only need to add either $u^{\text{out}}v^{\text{in}}$ or $v^{\text{out}}u^{\text{in}}$ to ${\cal G}$.
\end{proof}


\subsection{Proof of Theorem~\ref{thm:FPT_TCC_undir}}\label{app:FPT_TCC_undir}

\begin{proof}
Let $\mathcal{G}$ be a temporal graph and $k$ be a positive integer. We first prove items $1$ and $2$, namely, that there exist algorithms running in time:
\begin{itemize}
    \item $O(k^{k\cdot \tau}\cdot n)$ that decides whether there is a \TCC of size at least $k$; and
    \item $O(2^{k^\tau}\cdot n)$ that decides whether there is a \cTCC of size at least $k$. 
\end{itemize}

Denote by $F$ the graph obtained from the reachability digraph $\reach({\mathcal{G}})$ by removing all edges that are not symmetric and taking the underlying graph. Lemma~\ref{lemma:cliquecomponent} tells us that {\TCC} is equivalent to finding a clique of size at least $k$ in $F$, while {\cTCC} is equivalent to finding a clique $S$ of size at least $k$ such that ${\cal G}[S]$ is connected. 
Observe that if $\Delta(F)\le \Delta$, then the former can be solved by testing, for every $u\in V(F)$ and every $S\subseteq N_F(u)$ with $|S| = k-1$, whether $S\cup \{u\}$ is a clique in $F$; this takes time $O(\Delta^{k}\cdot k^2\cdot n)$. 
Now for the latter, we need to test for the existence of such sets of bigger sizes. This is because {\cTCC's} are not closed under inclusion. Nevertheless, since $\Delta(G)\le \Delta$ and testing whether ${\cal G}[S]$ is connected can be done in time $O(\lvert S\rvert\cdot \lvert E(G[S]))$, we can test for the existence of a ${\cTCC}$ in time $O(2^{\Delta}\cdot \Delta^3\cdot n)$ by searching all cliques of size at least $k$ in $N[u]$, for every $u\in V(G)$.
We finish the proof by bounding the value of $\Delta$.

Now, we show that $\Delta\le (k-1)^\tau$, which combined with the previous paragraph gives us the stated running time. For this, first notice that, for every $i\in[\tau]$, the vertex set of any connected component $C$ of $G_i$ is a clique in $F$ and a closed temporal connected set of $\mathcal{G}$.
This means that we can suppose that the size of any connected component of $G_i$ is at most $k-1$, for every $i\in [\tau]$, as otherwise we have a trivial yes-instance for both problems.
Now, since $C_i(u) = R_i(u)$ when $G$ is undirected, apply Lemma~\ref{lem:reachableset} to see that $\lvert {\cal R}_\tau(u)\rvert \le (k-1)^\tau$. 
Additionally, by definition we know that ${\cal R}_\tau(u)$ contains exactly the set of vertices reachable by $u$ in ${\cal G}$. Since $v\in N_F(u)$ if and only if $u$ reaches $v$ \emph{and} $v$ reaches $u$, it follows that $N_F(u)\subseteq {\cal R}_\tau(u)$. This finishes the proof of items 1 and 2.

Now we turn our attention to items 3 and 4, namely, algorithms running in time:
\begin{itemize}
    \item $O(k^{k^2}\cdot n)$ that decides whether there is a \TUCC of size at least $k$; and
    \item $O(2^{k^k}\cdot n)$ that decides whether there is a \cTUCC of size at least $k$.
\end{itemize}

Again by applying Lemma~\ref{lemma:cliquecomponent}, a similar argument as the one used for items $1$ and $2$  can be applied directly to the reachability graph $F = \reach({\cal G})$ to say that, if $\Delta(F)\le \Delta$, then {\TUCC} can be solved in time $O(\Delta^{k}\cdot k^2\cdot n)$, while {\cTUCC} can be solved in time $O(2^{\Delta}\cdot \Delta^3\cdot n)$. hence it remains to bound $\Delta$.

We prove first that $d_G(u)\le k-2$ for every $u\in V(G)$. This holds because, given any pair $v,w\in N(u)$, and any choice of values $i\in \lambda(uv)$ and $j\in \lambda(uw)$, either we have $i\le j$, in which case $(v,i,u,j,w)$ is a temporal path, or $i>j$, in which case $(w,j,u,i,v)$ is a temporal path. In other words, for every $v,w\in N(u)$, either $v$ reaches $w$ or $w$ reaches $v$, which implies that $N_G[u]$ is a clique in $F$, for every $u\in V(G)$. Hence if $d_G(u)\ge k-1$, we are in a trivial yes instance. To finish, just observe that any temporal path forms a {\cTUCC}, which is always contained in a {\TUCC}. Therefore, we can suppose that any vertex reachable from $u$ is reached by a temporal path containing at most $k-1$ edges. Because $d_G(v)\le k-2$ for every $v\in V(G)$, we get that $d_F(u) = \lvert {\cal R}_\tau(u)\rvert\le (k-2)^{k-1}$ and the result follows.
\end{proof}


\section{Checking Connectivity: proof of Theorem~\ref{thm:lowerbound_algorithm}}
\label{app:lower}

This section is focused on Question~\ref{quest:lower}, which is open for all definitions of components for both the strict and the non-strict models. We answer to the question providing the conditional lower bound in Theorem~\ref{thm:lowerbound_algorithm}, which holds for both models, where the notation $\Tilde{O}(\cdot)$ ignores poly-logarithmic factors.

We apply the technique used for instance in~\cite{BorassiCH16,puatracscu2010possibility,williams2010subcubic} to prove lower bounds for polynomial problems, falling within the fine-grained complexity framework.
%
We use \emph{quasilinear Karp reduction}, i.e. Karp reductions running in quasilinear time, whose formal definition is given in~\cite{BorassiCH16}. In the following we will use $\Tilde{O}(\cdot)$ to neglect poly-logarithmic factors. 

The key idea is to reduce a starting problem that is known not to be solvable in subquadratic time to our problem using such kind of reduction. 
This seed problem is the following formulation of the $k$-$\SAT^*$ problem. Let $\phi$ be a CNF formula on variables ${\cal X} = \{x_1,\ldots,x_n\}$ and ${\cal Y} = \{y_1,\ldots, y_n\}$, with $m$ clauses of size at most $k$. Let  $X$ denote the set of all $2^n$ possible truth assignments for ${\cal X}$, and similarly let $Y$ denote the set of all $2^n$ possible truth assignments for ${\cal Y}$. In the $k$-$\SAT^*$ problem, given $I = (\phi, X, Y)$, the goal is to decide if $\phi$ is satisfiable. The main difference with relation to the classical $k$-$\SAT$ problem is the size of the input, which is $|I|=O(2^n)$.

\begin{remark}[\cite{BorassiCH16}]\label{rem:SETH}
$k$-$\SAT^*$ with input $I$ cannot be solved in time $O(|I|^{2-\epsilon})$ for some $\epsilon$, unless \textsc{SETH} fails.
\end{remark}

By presenting a quasilinear Karp reduction from $k$-$\SAT^*$, and applying Remark~\ref{rem:SETH}, we obtain that, unless SETH fails, there is no subquadratic algorithm that decides if a given temporal graph is temporal (unilaterally) connected. 

In this section, it is helpful to formally define the following two problems.

\myprob{\textsc{Temporal Connected}}{A temporal graph ${\cal G}$.}{Is ${\cal G}$ temporal connected?}

\myprob{\textsc{Temporal Unilaterally  Connected}}{A temporal graph ${\cal G}$.}{Is ${\cal G}$ temporal unilaterally  connected?}

For both problems, given an instance $I = (\phi, X, Y)$ of $k$-$\SAT^*$, we construct a temporal graph ${\cal G} = (G,\lambda)$ such that ${\cal G}$ is not temporal (unilaterally) connected if and only if $\phi$ has a satisfying assignment. As in the obtained temporal graph all non-strict paths are also strict, the result holds on both models.

\begin{figure}[t]
\centering
\begin{center}
\begin{tikzpicture}[vertex/.style={circle, minimum size=0.2cm, draw, inner sep=1pt, fill=white}, arrow/.style={-{Stealth[scale=1.5]}, shorten >= 2pt}]

  \pgfsetlinewidth{1pt}
  \pgfdeclarelayer{bg}    
   \pgfsetlayers{bg,main}  

  \tikzset{snapshot/.style ={draw=black!50, rounded corners, dashed, minimum height=8mm, minimum width=2.3cm} }
  \tikzset{subgraph/.style ={draw=black!50, circle, x radius=.8cm, y radius=1.7cm, draw, minimum width=1cm, yscale=2, fill=white} }

    \node[vertex] (s) at (-2,0) {$s$};
    \node [subgraph] (X) at (0,0) {$X$};
    \node [subgraph] (C) at (3,0) {$C$};
    \node [subgraph] (Y) at (6,0) {$Y$};

    \begin{pgfonlayer}{bg}    
        \draw[arrow, thick] (X.-130) -- (s) node [midway,fill=white]{$6$};
        \draw[arrow, thick] (s) -- (X.130) node [midway,fill=white]{$7$};
        \draw[arrow, thick] (X.-130) -- (C.-130) node [midway,fill=white]{$8$};
        \draw[arrow, thick, shorten >= 1cm] (C.-160) to node [midway,fill=white]{$5$} (X.-160);
        \draw[arrow, red, thick] (X.130) -- (C.130) node [midway,fill=white]{$4$};
        \draw[arrow, thick] (C.-130) -- (Y.-130) node [midway,fill=white]{$3$};
        \draw[arrow, thick, shorten >= 1cm] (Y.-160) -- (C.-160) node [midway,fill=white]{$2$};
        \draw[arrow, red, thick] (C.130) -- (Y.130) node [midway,fill=white]{$5$};
        \node (hook) at (3, -2) {};
        \draw[thick] (s) to [out = -90, in = 180] (hook);
        \draw[arrow,thick] (hook) to [out = 0, in = -90] (Y);
        \node (name, fill = white) at (hook) {$1$};
    \end{pgfonlayer}
    
  \end{tikzpicture}

\end{center}
\vspace{-.5cm}
\caption{General structure of the constructed graph in the reduction for the temporal connectivity testing problem. A black edge denotes the existence of all possible edges. A red edge denotes the temporal edges whose existence is conditioned to the assignment not satisfying the clause. }
\label{fig:TCC_square}
\end{figure}
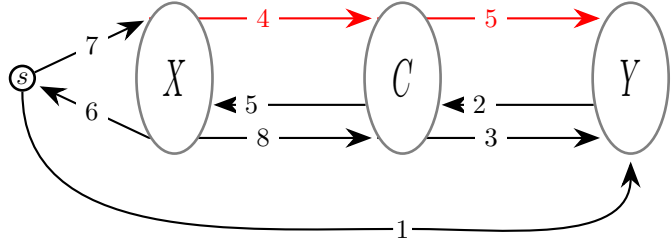

\begin{figure}[h]
\begin{center}
\begin{tikzpicture}[vertex/.style={circle, minimum size=0.2cm, draw, inner sep=1pt, fill=white}, arrow/.style={-{Stealth[scale=1.5]}, shorten >= 2pt}]

  \pgfsetlinewidth{1pt}
  \pgfdeclarelayer{bg}    
   \pgfsetlayers{bg,main}  

  \tikzset{snapshot/.style ={draw=black!50, rounded corners, dashed, minimum height=8mm, minimum width=2.3cm} }
  \tikzset{subgraph/.style ={draw=black!50, circle, x radius=.8cm, y radius=1.7cm, draw, minimum width=1cm, yscale=2, fill=white} }

    \node[vertex,draw=blue] (x1x2) at (0,2) {$TT$};
    \node[vertex,draw=blue] (x1nx2) at (0,1) {$TF$};
    \node[vertex,draw=blue] (nx1x2) at (0,0) {$FT$};
    \node[vertex,draw=blue] (nx1nx2) at (0,-1) {$FF$};

    \node[vertex,draw=magenta] (y1y2) at (6,2) {$TT$};
    \node[vertex,draw=magenta] (y1ny2) at (6,1) {$TF$};
    \node[vertex,draw=magenta] (ny1y2) at (6,0) {$FT$};
    \node[vertex,draw=magenta] (ny1ny2) at (6,-1) {$FF$};

    \node[vertex] (c1) at (3,2) {$c_1$};
    \node[vertex] (c2) at (3,0.5) {$c_2$};
    \node[vertex] (c3) at (3,-1) {$c_3$};
    
    \node[vertex] (s) at (-3,0.5) {$s$};
    
    \begin{pgfonlayer}{bg}    
        \draw[arrow, thick] (2.5,-2) -- (0.5,-2) node [midway,fill=white]{$5$};
        \draw[arrow, thick] (0.5,-2.5) -- (2.5,-2.5) node [midway,fill=white]{$8$};
        \draw[arrow, thick] (5.5,-2) -- (3.5,-2) node [midway,fill=white]{$2$};
        \draw[arrow, thick] (3.5,-2.5) -- (5.5,-2.5) node [midway,fill=white]{$3$};
        \draw[arrow, thick] (s) to [out=-90,in=-90] node [midway,fill=white]{$1$} (6,-1.5);
        \draw[arrow, thick] (s) to node [midway,fill=white]{$7$} (-0.5,2);
        \draw[arrow, thick] (-0.5,-1) to node [midway,fill=white]{$6$} (s);

        \draw[arrow, thick, blue] (x1x2) to  (c2);
        \draw[arrow, thick,blue] (x1nx2) to (c2);
        \draw[arrow, thick,blue] (x1nx2) to (c3);
        \draw[arrow, thick,blue] (nx1x2) to (c1);
        \draw[arrow, thick,blue] (nx1nx2) to (c3);

        \draw[arrow, thick, magenta] (c2) to (y1ny2);
        \draw[arrow, thick, magenta] (c1) to (ny1y2);
        \draw[arrow, thick, magenta] (c3) to (ny1y2);
        \draw[arrow, thick, magenta] (c1) to (ny1ny2);
    \end{pgfonlayer}
    
  \end{tikzpicture}

\end{center}
\caption{Graph in the reduction for the temporal connectivity testing problem, related to the formula $\phi = (x_1\vee \neg x_2 \vee y_1)\wedge (\neg x_1\vee \neg y_1\vee y_2)\wedge (x_2\vee y_1\vee \neg y_2)$. Blue nodes denote assignments of $x_1$ and $x_2$ (e.g., node $TT$ blue denotes the assignment $x_1 = True$ and $x_2 = True$), while magenta nodes denote assignments of $y_1$ and $y_2$. Again, black edges denote existence of all possible edges. We put them outside the vertices in order to make the figure clean. Blue edges are active in time 4 and magenta edges, in time 5.}
\label{fig:TCC_square_example}
\end{figure}
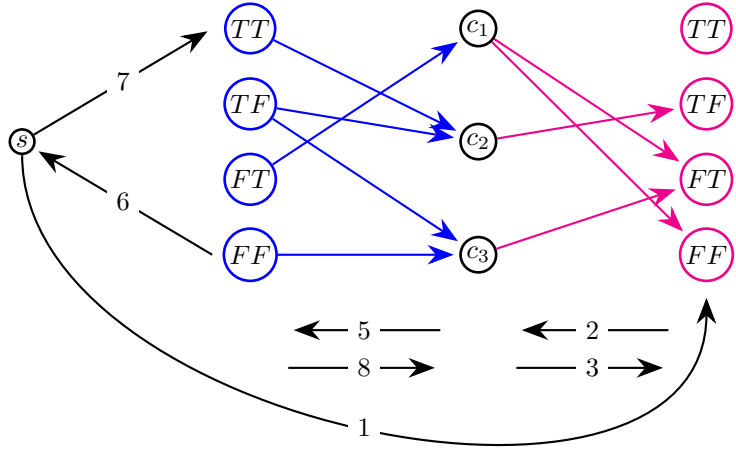

We first present a reduction from $k$-$\SAT^*$ to the complement of \textsc{Temporal Connected}.
See Figure~\ref{fig:TCC_square} to follow the construction. Also, let $C = \{c_1,\ldots,c_m\}$ be the set of clauses in $\phi$. Let $V(G) = X\cup C\cup Y\cup \{s\}$. Add all edges from $s$ to $X$ and let them be active in time $7$; all edges between $X$ and $s$ and let them be active in time $6$; all edges between $X$ and $C$, active in time $8$; all edges between $C$ and $X$, active in time $5$; all edges between $C$ and $Y$, active in time $3$; all edges between $Y$ and $C$, active in time $2$; and all edges between $s$ and $Y$, active in time $1$. Finally, for each pair $f\in X$ and $c_i\in C$, we add edge $fc_i$, active in time~4, if and only if $f$ does \emph{not} satisfy $c_i$. Similarly, for each pair $f\in Y$ and $c_i\in C$, we add edge $c_if$, active in time~5, if and only if $f$ does \emph{not} satisfy $c_i$. See Figure~\ref{fig:TCC_square_example} for an example. 

We now argue that this is a quasilinear Karp reduction. 
The reduction is quasilinear since $\tau = 8$, $\lvert V(G)\rvert = 2^n+2^n+m+1 = O(\lvert I \rvert)$, and $\lvert E(G) \rvert = 2^nm+2^nm+2^n+2^n = \Tilde{O}(\lvert I\rvert)$, since $m$ in Remark~\ref{rem:SETH} can be assumed to be $polylog(\lvert I\rvert)$~\cite{BorassiCH16}. It remains to prove correctness. Before we do that, we first argue that the reachability graph of ${\cal G}$ always contains $uv$ for every $u$ and $v$ such that either $u\notin X$ or $v\notin Y$. For this, we analyse all cases below:

\begin{itemize}
    \item $s$ reaches $f$ for every $f\in X\cup Y$ by direct edges;
    \item $s$ reaches $c_i$ for every $c_i\in C$ through a path $(s,1,f,2,c_i)$ for any $f\in Y$;
    \item $s$ is reachable by $f$ for every $f\in X$ by direct edges;
    \item $s$ is reachable by every $c_i\in C$ through a path $(c_i,5,f,6,s)$ for any $f\in X$;
    \item $s$ is reachable by every $f\in Y$ through a path $(f,2,c_i,5,f',6,s)$ for any $c_i\in C_i$ and any $f'\in X$;
    \item Every $f\in X$ reaches every $c_i\in C$ and is reached by it through direct edges;
    \item Every $f\in X$ is reachable by every $f'\in Y$ through a path $(f',2,c_i,5,f)$ for any $c_i\in C_i$;
    \item Every $c_i\in C$ reaches every $f\in Y$ and is reached by it through direct edges;
    \item Every $f\in X$ reaches every $f'\in X$ through the path $(f,6,s,7,f')$;
    \item Every $c_i\in C$ reaches every $c_j\in C$ through the path $(c_i,5,f,8,c_j)$ for any $f\in X$;
    \item Every $f\in Y$ reaches every $f'\in Y$ through the path $(f,2,c_i,3,f')$ for any $c_i\in C$.
\end{itemize}

Now we prove that $\phi$ is satisfiable if and only if there exists $f_X\in X$ and $f_Y\in Y$ such that $f_X$ does \emph{not} reach $f_Y$ (i.e., ${\cal G}$ is \emph{not} temporal connected). First, suppose that $\phi$ is satisfiable and consider a satisfying assignment $f$ of $\phi$. Then let $f_X\in X$ be equal to $f$ restricted to $\{x_1,\ldots, x_n\}$ and $f_Y\in Y$ be equal to $f$ restricted to $\{y_1,\ldots, y_n\}$. Observe that for every $c_i\in C$, either $f_X$ satisfies $c_i$, and hence $(f_Xc_i,4)\notin E^T({\cal G})$, or $f_Y$ satisfies $c_i$, and hence $(c_if_Y,5)\notin E^T({\cal G})$. Observe also that the only possible temporal paths between $f_X$ and $f_Y$ are of the type $(f_X,4,c_i,5,f_Y)$ for some $c_i\in C$. It thus follows that $f_X$ does not reach $f_Y$. Now suppose that $f_X$ does not reach $f_Y$ for some pair $f_X\in X$ and $f_Y\in Y$. This must be because for every $c_i\in C$, either $(f_Xc_i,4)\notin E^T({\cal G})$, and hence $f_X$ satisfies $c_i$, or $(c_if_Y,5)\notin E^T({\cal G})$, and hence $f_Y$ satisfies $c_i$. Therefore $f_X\cup f_Y$ is a satisfying assignment for $\phi$.

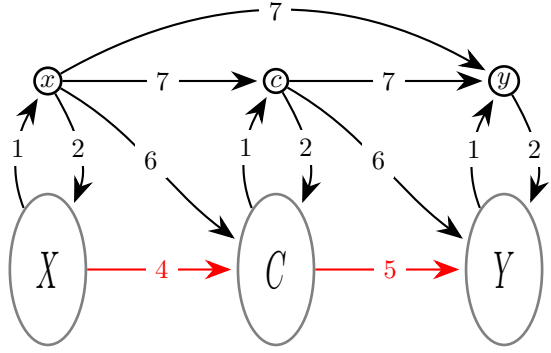
\begin{figure}[t]
\centering
\begin{tikzpicture}[vertex/.style={circle, minimum size=0.2cm, draw, inner sep=1pt, fill=white}, arrow/.style={-{Stealth[scale=1.5]}, shorten >= 2pt}]

  \pgfsetlinewidth{1pt}
  \pgfdeclarelayer{bg}    
   \pgfsetlayers{bg,main}  

  \tikzset{snapshot/.style ={draw=black!50, rounded corners, dashed, minimum height=8mm, minimum width=2.3cm} }
  \tikzset{subgraph/.style ={draw=black!50, circle, x radius=.8cm, y radius=1.7cm, draw, minimum width=1cm, yscale=2, fill=white} }

    \node [subgraph] (X) at (0,0) {$X$};
    \node[vertex] (x) at (0,2.5) {$x$};
    \node [subgraph] (C) at (3,0) {$C$};
    \node[vertex] (c) at (3,2.5) {$c$};
    \node [subgraph] (Y) at (6,0) {$Y$};
    \node[vertex] (y) at (6,2.5) {$y$};

    \begin{pgfonlayer}{bg}    
        \draw[arrow, thick] (x) to [out=-60,in=60] node [midway,fill=white]{$2$} (X);
        \draw[arrow, thick] (x) to node [midway,fill=white]{$7$} (c);
        \draw[arrow, thick] (x) to [out=30,in=150] node [midway,fill=white]{$7$} (y);
        \draw[arrow, thick] (c) to node [midway,fill=white]{$7$} (y);

        \draw[arrow, thick] (X) to [out=120,in=-120] node [midway,fill=white]{$1$} (x);
        \draw[arrow, red, thick] (X) -- (C) node [midway,fill=white]{$4$};
        \draw[arrow, thick] (x) to [out=-30,in=150] node [midway,fill=white]{$6$} (C);
        
        \draw[arrow, thick] (c) to [out=-60,in=60] node [midway,fill=white]{$2$} (C);
        \draw[arrow, thick] (C) to [out=120,in=-120] node [midway,fill=white]{$1$} (c);
        \draw[arrow, red, thick] (C) -- (Y) node [midway,fill=white]{$5$};
        \draw[arrow, thick] (c) to [out=-30,in=150] node [midway,fill=white]{$6$} (Y);

        \draw[arrow, thick] (y) to [out=-60,in=60] node [midway,fill=white]{$2$} (Y);
        \draw[arrow, thick] (Y) to [out=120,in=-120] node [midway,fill=white]{$1$} (y);
        
    \end{pgfonlayer}
  \end{tikzpicture} 
\caption{General structure of the constructed graph in the reduction for the unilateral temporal connectivity testing problem. A black edge denotes the existence of all possible edges. A red edge denotes the temporal edges whose existence is conditioned to the assignment not satisfying the clause.}
\label{fig:TUCC_square}
\end{figure}

Consider now the complement of \textsc{Temporally Unilaterally Connected}. We make a similar reduction. Observe Figure~\ref{fig:TUCC_square} to follow the construction. 
Let $V(G) = X\cup C\cup Y\cup \{x,y,c\}$. For each $z\in \{x,y,c\}$, add all edges from $Z$ to $z$ and let them be active in time~$1$, and all edges from $z$ to $Z$ and let them be active in time~$2$. Add also all edges from $x$ to $C$ and from $c$ to $Y$, active in time~$6$. Finally, add $\{xc,xy,cy\}$ active in time $7$ and, for each pair $f\in X$ and $c_i\in C$, we add edge $fc_i$, active in time~$4$, if and only if $f$ does \emph{not} satisfy $c_i$. Similarly, for each pair $f\in Y$ and $c_i\in C$, we add edge $c_if$, active in time~$5$, if and only if $f$ does \emph{not} satisfy $c_i$. 

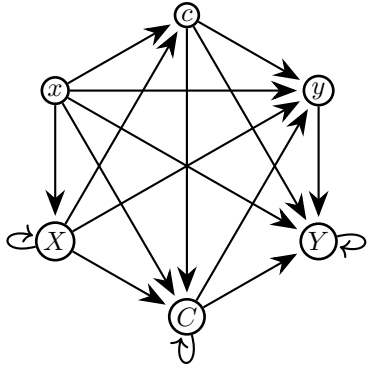
\begin{figure}
\begin{center}
\begin{tikzpicture}[vertex/.style={circle, minimum size=0.2cm, draw, inner sep=1pt, fill=white}, arrow/.style={-{Stealth[scale=1.5]}, shorten >= 2pt}]

  \pgfsetlinewidth{1pt}
  \pgfdeclarelayer{bg}    
   \pgfsetlayers{bg,main}  

  \tikzset{snapshot/.style ={draw=black!50, rounded corners, dashed, minimum height=8mm, minimum width=2.3cm} }
  \tikzset{subgraph/.style ={draw=black!50, circle, x radius=.8cm, y radius=1.7cm, draw, minimum width=1cm, yscale=2, fill=white} }

    \node [vertex] (X) at (-150:2) {$X$};
    \node[vertex] (x) at (150:2) {$x$};
    \node [vertex] (C) at (-90:2) {$C$};
    \node[vertex] (c) at (90:2) {$c$};
    \node [vertex] (Y) at (-30:2) {$Y$};
    \node[vertex] (y) at (30:2) {$y$};

    \begin{pgfonlayer}{bg}    
        
        \draw[arrow,thick] (x)--(X);
        \draw[arrow,thick] (x)--(y);
        \draw[arrow,thick] (x)--(c);
        \draw[arrow,thick] (x)--(Y);
        \draw[arrow,thick] (x)--(C);

        \draw[arrow,thick] (X)--(c);
        \draw[arrow,thick] (X)--(y);
        \draw[arrow,thick] (X)--(C);
        \path[arrow,thick] (X) edge  [loop left] node {} ();

        \draw[arrow,thick] (c)--(Y);
        \draw[arrow,thick] (c)--(y);
        \draw[arrow,thick] (c)--(C);

        \draw[arrow,thick] (C)--(y);
        \draw[arrow,thick] (C)--(Y);
        \path[arrow,thick] (C) edge  [loop below] node {} ();

        \draw[arrow,thick] (y)--(Y);
        \path[arrow,thick] (Y) edge  [loop right] node {} ();

    \end{pgfonlayer}
    
  \end{tikzpicture}

\end{center}
\caption{Subgraph of $\reach({\cal G})$ shown above. An edge between two sets/vertices means that all such edges exist in $\reach({\cal G})$. Note that more edges might exist but that this is enough to prove that we miss only an edge between $X$ and $Y$ (which means all edges between these sets).}
\label{fig:TUCC_square_reachability}
\end{figure}

Let ${\cal G} = (G,\lambda)$ be the constructed temporal graph. Similarly as before, we get $\tau = 7$, $\lvert V(G)\rvert = 2^{n+1}+m+3 = O(\lvert I\rvert)$, and $\lvert E(G)\rvert = 2^{n+1}m + 2^{n+1} + m + 3 = \Tilde{O}(\lvert I\rvert)$. It remains to prove that $\phi$ is not satisfiable if and only if ${\cal G}$ is unilaterally temporal connected. As before, we first prove that the only missing pairs are of the type $f_Xf_Y$ with $f_X\in X$ and $f_Y\in Y$. Recall that we only need at least one of the edges $uv$ or $vu$ for every pair $u,v\in V(G)$. Below we analyse only the necessary edges. See Figure~\ref{fig:TUCC_square_reachability} to follow the proof.

\begin{itemize}
    \item $x$ reaches every $u\in V(G)$ directly through edges. This implies edges $xu$ in $\reach({\cal G})$ for every $u\in V(G)$;
    \item Every $f_X\in X$ reaches every $z\in \{c,y\}$ through the path $(f_X,1,x,7,z)$. This implies edges $f_Xc$ and $f_Xy$ in $\reach({\cal G})$;
    \item Every $f_X\in X$ reaches every $c_i\in C$ through the path $(f_X, 1, x, 6, c_i)$. This implies edges $f_Xc$ in $\reach({\cal G})$;
    \item Every $f_X\in X$ reaches every $f'_X\in X$, $f'_X\neq f_X$, through the path $(f_X, 1, x, 2, f'_X)$;
    \item $c$ reaches every $u\in C\cup Y\cup \{y\}$ directly through edges. This implies edges $cu$ in $\reach({\cal G})$;
    \item Every $c_i\in C$ reaches $y$ through the path $(c_i,1,c,7,y)$. This implies edges $c_iy$ in $\reach({\cal G})$;
    \item Every $c_i\in C$ reaches every $f_Y\in Y$ through the path $(c_i, 1, c, 6, f_Y)$. This implies edges $c_if_Y$ in $\reach({\cal G})$;
    \item Every $c_i\in C$ reaches every $c_j\in C$, $i\neq j$, through the path $(c_i, 1, c, 2, c_j)$. This implies edges $c_ic_j$ in $\reach({\cal G})$;
    \item $y$ reaches every $f_Y\in Y$ directly  through an edge. This implies edges $yf_Y$ in $\reach({\cal G})$;
    \item Every $f_Y\in Y$ reaches every $f'_Y\in Y$, $f'_Y\neq f_Y$, through the path $(f_Y, 1, y, 2, f'_Y)$. This implies edges $f_Yf'_Y$ in $\reach({\cal G})$.
\end{itemize}

The proof of correctness is analogous to the previous one. It relies on the fact that the only possible temporal paths between the sets $X$ and $Y$ are of the type $(f_X,4,c_i,5,f_Y)$, for some $c_i\in C$. Indeed, no path from $Y$ to $X$ may exist, since each edge that leaves $Y$ goes to $y$, and there are no edges from $y$ to $V(G)\setminus Y$. Additionally, the only temporal edges leaving $X$ that do not go directly to $C$ are of type $(f_Xx,1)$. If a temporal path starting with such temporal edge uses an edge at time 7, then it cannot arrive to $Y$, as all edges arriving in $Y$ have smaller timestamps. Hence such a path must use an edge $(xc_i,6)$ for some $c_i\in C$, and again we get stuck as all edges leaving $C$ occur before time~$6$. 
Now, since the path $(f_X,4,c_i,5,f_Y)$ exists if and only if both $f_X$ and $f_Y$ do not satisfy $c_i$, the results follows. That is, if $\phi$ is satisfiable, then this fails for some pair, and vice-versa. 

Now observe that we have made reductions from $k$-$\SAT^*$ to the \emph{complements} of our problems. However, since a subquadratic algorithm that solves the complement $\overline{\Pi}$ of a problem $\Pi$, also solves $\Pi$ (indeed $I$ is a positive instance of $\Pi$ if and only if $I$ is a negative instance of $\overline{\Pi}$), we get that Theorem~\ref{thm:lowerbound_algorithm} follows.


\section{Checking Maximality: proof of Theorem~\ref{thm:maximality}}\label{app:maximality}

\begin{proof}
It remains to prove that, in the non-strict model, given a (directed) temporal graph $\mathcal{G}$ and a subset $Y\subseteq V(\mathcal{G})$, deciding whether $Y$ is a \cTCC (\cTUCC) is $\NP$-complete. The strict case is already treated in the main text. As before, we make a reduction from the problem of, given a graph $G$ and $X \subseteq V(G)$, deciding whether $X$ is a $2$-club. 

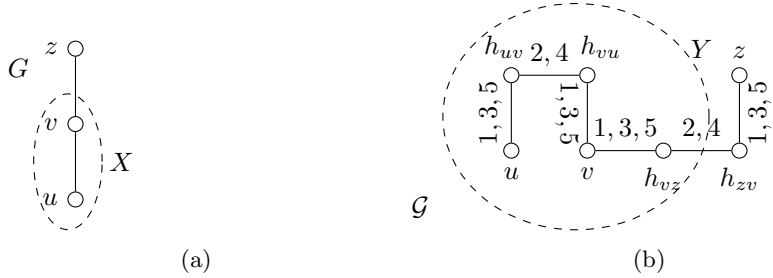
\begin{figure}[t]
    \centering
    \begin{subfigure}[b]{0.4\textwidth}
    {\begin{tikzpicture}[scale=1]

\node (G) at (-.75,1.75) {$G$};
\node[dashed,draw, ellipse, minimum width=.9cm, minimum height=1.8 cm] at (-.1,.5) {};
\node[label=right:$X$] at (0.2,.5) {};

\vertex[label=left:$u$] (u) at (0,0) {};
\vertex[label=left:$v$] (v) at (0,1) {};
\vertex[label=left:$z$] (z) at (0,2) {};

\draw (u) -- (v) -- (z);

\end{tikzpicture}}
    \caption{}
    \end{subfigure}
    \begin{subfigure}[b]{0.5\textwidth}
    {\begin{tikzpicture}[scale=1]
 \node[dashed,draw,ellipse, minimum width=3.5cm,minimum height=3cm] at (.85,.45){};
 \node (G) at (-1.2,-.75) {$\mathcal{G}$};
\node (Y) at (2.5,1.35) {$Y$}{};
\vertex[label=below:$u$] (u) at (0,0) {};
\vertex[label=above:$\!\!\!h_{uv}$] (huv) at (0,1) {};
\vertex[label=above:$~~~h_{vu}$] (uvh) at (1,1) {};
\vertex[label=below:$v$] (v) at (1,0) {};
\vertex[label=below:$h_{vz}$] (hvz) at (2,0) {};
\vertex[label=below:$h_{zv}$] (vzh) at (3,0) {};
\vertex[label=above:$z$] (z) at (3,1) {};

\draw (u) to node[scale=1,anchor=south,sloped, above]{$1,3,5$}  (huv);
\draw (huv) to node[scale=1,sloped, above] {$2,4$} (uvh);
\draw (uvh) to node[scale=1,anchor=north,sloped, below]{$1,3,5$} (v);
\draw (v) to node[scale=1,sloped, above] {$1,3,5$} (hvz);
\draw (hvz) to node[scale=1,sloped, above] {$2,4$} (vzh);
\draw (vzh) to node[scale=1,anchor=north,sloped, below]{$1,3,5$} (z);

\end{tikzpicture}}
    \caption{}
    \end{subfigure}
    \caption{Construction in the proof of Theorem~\ref{thm:maximality}.}
    \label{fig_app:thm7}
\end{figure}

Observe Figure~\ref{fig_app:thm7}. We obtain $\mathcal{G}$ from $G$ by subdividing each edge $uv \in E(G)$ twice, creating vertices $h_{uv}$ and $h_{vu}$, with $\lambda(uh_{uv}) = \lambda(vh_{vu}) = \{1,3,5\}$, and $\lambda(h_{uv}h_{vu}) = \{2,4\}$. 
Observe that the first vertex in the subscript of $h_{xy}$ tells us which between $x$ and $y$ is adjacent to $h_{xy}$. 
Denote by $H$ the set $\{h_{uv},h_{vu}\mid uv\in E(G)\}$. We now prove that $X$ is a maximal $2$-club in $G$ if and only if $Y = X\cup N_H(X)$ is a \cTCC in $\mathcal{G}$. In fact, we prove that:
\begin{enumerate}
    \item\label{1} If $X\subseteq V(G)$ is such that $G[X]$ has diameter at most~2, then $Y = X\cup N_H(X)$ is a closed connected set in $\mathcal{G}$; and
    \item\label{2} If $Y\subseteq V(\mathcal{G})$ is a closed connected set, then $X = Y\cap V(G)$ is such that $G[X]$ has diameter at most~2.
\end{enumerate}
 
 We argue that indeed \ref{1} and \ref{2} above imply what we want, i.e., that $X$ is a maximal $2$-club in $G$ if and only if $Y = X\cup N_H(X)$ is a \cTCC in $\mathcal{G}$. Observe that, supposing that~\ref{1} and~\ref{2} hold, if $X$ is a maximal 2-club, then $Y = X\cup N_H(X)$ must be a \cTCC. Indeed, if $Y\subset Y'$ and $Y'$ is a closed connected set (i.e., $Y$ is not maximal), then by~\ref{2} we get that $X' = Y'\cap V(G)$ has diameter~2. Since $X'$ contains $X$, this contradicts the choice of $X$. Conversely, if $Y$ is a \cTCC, then $X$ must be a maximal 2-club, as otherwise we could apply~\ref{1} to get a closed connected set strictly containing $Y$. 

We first prove~\ref{1}. So, consider $X\subseteq V(G)$ such that $G[X]$ has diameter at most~2, and define $Y$ as above. 
Let $u,v\in Y\cap V(G)$. If $uv\in E(G)$, then $(u,1,h_{uv},2,h_{vu},3,v)$ and $(v,1,h_{vu},2,h_{uv},3,u)$ witness that $u$ reaches $v$ and $v$ reaches $u$ in $\mathcal{G}[Y]$. 
And if $uv\notin E(G)$, then, since $G[X]$ has diameter~2, let  $w\in N(u)\cap N(v)$ in $G$. 
We get the paths: $(u,1,h_{uw},2,h_{wu},3,w,3,h_{wv},4,h_{vw},5,v)$ and $(v,1,h_{vw},2,h_{wv},3,w,3,h_{wu},$ $4,h_{uw},5,u)$. Therefore, $u$ reaches $v$ and $v$ reaches $u$ in $\mathcal{G}[Y]$. 
Now, consider $u\in X\cap Y$ and $h\in H\cap Y$. Let $v\in X$ be such that $h\in N(v)$ (observe that $v$ is uniquely defined). If $v=u$, then there is nothing to prove, so suppose otherwise. 
Because $X$ has diameter at most~2, there exists a $u,v$-path $P$ in $G[X]$ of length at most~2, say $(u,w,v)$, with possibly $w=v$. Then either $(u,1,h_{uw},2,h_{wu},3,w,3,h)$ is a temporal $u,h$-path, in case $w=v$, or $(u,1,h_{uw},2,h_{wu},3,w,3,h_{wv},4,h_{vw},5,v,5,h)$ is a temporal $u,h$-path, in case $w\neq v$. One can check that the symmetric path between $h$ and $u$ ensures that also $h$ reaches $u$ in $\mathcal{G}[Y]$. Now, let $h,h'\in H$, and let $u\in N(h)\cap X$ and $v\in N(h')\cap X$. One can observe that a similar argument can be applied, by possibly starting the previous path with $(h,1,u)$, in case $h = h_{ux}$ for some $x$ not within the $u,v$-path $P$ taken in $G[X]$.

Now, assume that $Y$ is a closed connected set of $\mathcal{G}$, and consider $X=Y\cap V(G)$. We want to show that $G[X]$ has diameter at most~2. 
Suppose by contradiction that $u$ and $v$ are a distance~$3$ in $G[X]$. 
Observe that, since each $h\in H$ has degree exactly~2 in $\mathcal{G}$, we get that every temporal path in $\mathcal{G}$ is related to exactly one path in $G$, and vice-versa. One can then verify that $u$ and $v$ cannot reach each other in $\mathcal{G}[Y]$ as the traversal of any edge in $G$ is related to the traversal of a strictly increasing path of length~3 in $\mathcal{G}[Y]$. In other words, no $u,v$-path in $\mathcal{G}$ is a valid temporal path, as the lifetime of $\mathcal{G}$ is~5.

Finally, we prove that every closed unilaterally connected set is also a closed connected set. Since the reverse trivially holds, we get that it is also $\NP$-complete to decide whether $Y\subseteq V(\mathcal{G})$ is a \cTUCC. So, consider $Y\subseteq V(\mathcal{G})$ a closed unilaterally connected set, and suppose that $x,y\in Y$ are such that $x$ reaches $y$ in $\mathcal{G}[Y]$. Let $(x = x_1,t_1,x_2,\ldots,x_q,t_q,x_{q+1}=y)$ be a temporal $x,y$-path in $\mathcal{G}[Y]$.  By a case analysis, one can check that there exist $t'_1,\ldots,t'_q$ for which $(y = x_{q+1},t'_1,x_{q-1},\ldots,x_2,t'_q,$ $x_{1}=x)$ is a valid temporal $y,x$-path.
\end{proof}
\end{document}